\documentclass[12pt]{article}
\usepackage{latexsym}
\usepackage{epsfig}
\usepackage{verbatim}
\usepackage{amssymb}
\usepackage{amsmath}

\title{Multifractional, multistable, and other processes
with prescribed local form}
\author{K.J. Falconer \\
\small{{\em Mathematical Institute,
University of St~Andrews, North Haugh, St~Andrews,}} \\
\small{{\em Fife, KY16~9SS, Scotland }} \\
\small{ and } \\
J. L\'{e}vy V\'{e}hel \\
\small{{\em Projet Fractales, INRIA Rocquencourt, 78153 Le Chesnay
Cedex, France}}}
\date{}

\newtheorem{theo}{Theorem}

\newtheorem{prop}[theo]{Proposition}

\newtheorem{lem}[theo]{Lemma}
\newtheorem{cor}[theo]{Corollary}

\newcommand{\bbbr}{\mathbb R}
\newcommand{\1}{{\bf 1}}

\newcommand\var{\mbox{\rm var}}

\newcommand\E{\mbox{\sf E}}

\newcommand\e{\eta}

\newcommand\cd{\stackrel{{\rm d}}{\rightarrow}}
\newcommand\cp{\stackrel{{\rm p}}{\rightarrow}}

\newcommand\fdd{\stackrel{{\rm fdd}}{\rightarrow}}

\renewcommand\P{{\sf P}}

\newcommand{\one}{\ifmmode {\sf 1}\hspace{-.26em}{\sf
l}\hspace{-.35em}{\sf \_} \else ${\sf 1}\hspace{-.26em}{\sf
l}\hspace{-.35em}{\sf \_}$ \fi}

\newcommand\s{\Sigma}
\renewcommand\th{\theta}
\newcommand\X{{\sf X}}
\newcommand\Y{{\sf Y}}

\renewcommand{\Box}{\mbox{\rule{1ex}{1ex}}}

\begin{document}
\maketitle
\begin{abstract}
\noindent We present a general method for constructing stochastic processes with
prescribed local form. Such processes include
variable amplitude multifractional Brownian motion,
multifractional $\alpha$-stable processes, and multistable processes,
that is processes that are locally $\alpha(t)$-stable but where 
the stability index $\alpha(t)$ varies with $t$.  In particular we
construct multifractional multistable processes, where both the local 
self-similarity and stability indices vary.     
    
\end{abstract}


\section{Introduction}
\setcounter{equation}{0}
\setcounter{theo}{0}
\medskip

In this paper we present a general framework for constructing
stochastic processes with prescribed local forms.

Stochastic processes where the local H\"{o}lder regularity 
varies with a parameter $t$ (usually time) are important 
both in theory and in practical applications.   The best known example
is multifractional 
Brownian motion (mBm), where the Hurst index $h$ of fractional
Brownian motion is replaced by a functional 
parameter $h(t)$, permitting the H\"{o}lder exponent to 
vary in a prescribed manner.  This allows local regularity 
and long range dependence to be decoupled 
to give sample paths that are both highly 
irregular and highly correlated, a useful feature
for terrain or TCP traffic modeling.

For modelling financial or medical data
another feature is often important, namely
the presence of jumps.  
Stable non-Gaussian processes give good models for 
data containing discontinuities, with the stability 
index $\alpha$ controlling the distribution of jumps.  
Recently, multifractional stable processes, 
generalising mBm, were introduced to provide
jump processes with varying local regularity.  
However, a further step is needed for situations where 
both local regularity and jump intensity vary with time, 
for example to model financial data or epileptic episodes in EEG,
where for some periods there may be
only small jumps and at other instants very large ones.  
Our method may be used to construct processes where both $h$ and $\alpha$ 
vary in a prescribed way: thus there are two parameters which might 
correspond to distinct aspects of financial risk, 
to different sources of irregularity leading 
to the onset of epilepsy, or to textured images
where both H\"{o}lder regularity and 
the distribution of discontinuities  varies.

It is natural to construct processes $Y =\{Y(t): t \in \bbbr\}$ that
have an identifiable local form near each $u$, that is
 where there is a limiting process
\begin{equation}
\lim_{r \to 0}\frac{Y(u+rt) -Y(u)}{r^{h}} = Y_{u}'(t)
\label{locform1}
\end{equation}
which may vary with $u$. If this limit exists  as a non-trivial process 
we will say that $Y$ is $h$-{\it localisable} at $u$
and call the process $Y_{u}'=\{Y_{u}'(t): t\in \bbbr\}$ the {\it local form} of $Y$ at $u$.
The limit (\ref{locform1}) may
be taken in several ways: of particular interest are  
convergence in finite dimensional distributions, and 
convergence in distribtion;
in the latter case we term the process {\it strongly $h$-localisable}.
We will be especially concerned with $h$-localisable processes with $0<h<1$ 
which are usually of a fractal nature.

The most familiar example is multifractional Brownian motion $Y$
 which resembles index-$h(u)$ fractional
Brownian motion close to time $u$ but where $h(u)$ varies, that is  
\begin{equation}
\lim_{r \to 0}\frac{Y(u+rt) -Y(u)}{r^{h}} = B_{h(u)}(t) \label{exfbm}
\end{equation}
where $B_{h}$ is index-$h$ fractional Brownian motion, see
\cite{AA,AL2,BJR,EH,PL}. 
Generalising this, mulitfractional 
$\alpha$-stable processes have been
constructed with local form $h(u)$-self-similar linear
$\alpha$-stable motions \cite{ST1,ST2}.

It is clear  from (\ref{locform1}) that the $h$-local form $Y_{u}'$ at
$u$, if it exists, 
must itself be $h$-self-similar, that is $Y_{u}'(rt)= r^{h}Y_{u}'(t)$ 
for $r>0$.
However, much more is true: under quite general conditions $Y_{u}'$ must 
be self-similar with stationary increments (sssi) at almost all $u$ at
which it is strongly localisable, that is  
$r^{-h}(Y_{u}'(u+rt) -Y_{u}'(u)) = Y_{u}'(t)$ for
all $u$ and $r>0$, see \cite{Fal5,Fal6}.  Thus if we 
wish to construct processes with given local forms, the local forms
should themselves be sssi.  Whilst this is a strong requirement, many
classes of sssi processes are known, including fractional Brownian
motion, linear fractional stable motion and $\alpha$-stable L\'{e}vy 
motion, see \cite{EM,ST}.

Our general construction will allow known localisable
processes $X(\cdot,v)=\{X(t,v): t\in \bbbr\}$ for a range of $v$ 
to be pieced together to yield a localisable
`diagonal' process 
$Y=\{X(t,t): t\in \bbbr\}$ with
local form depending on $t$.  We will obtain conditions for the
transference of the local properties of $X(\cdot,v)$ to $Y$. 
The basic setting is akin to that
adopted in \cite{AL2,ST1}.  Thus we seek a random field
$\{X(t,v):(t,v) \in \bbbr^{2}\}$ such that
for each $v$ the local form $X_{v}'(\cdot,v)$ of $X(\cdot,v)$ at $v$
is the desired local 
form $Y_{v}'$ of $Y$ at $v$.  Typically, for each $v$ the
process $\{X(t,v): t\in \bbbr\}$ will be one where the local form can be readily
identified, such as an sssi process. Clearly the interplay of 
$X(\cdot,v)$ for $v$ in a neighbourhood of $u$ will be crucial to the local
behaviour of $Y$ near $u$.  
Thus the random field is set up as
an integral or sum of functions that depend on $t$ and $v$ with
respect to a single underlying random measure or process to
provide the necessary correlations.  In Section 4 we derive
general criteria that guarantee the transference of
localisability from the $X(\cdot,v)$ to $Y=\{X(t,t): t \in \bbbr\}$; Section 5
addresses this for strong localisability.

We illustrate the general method with several specific classes of
processes.
The method permits easy constructions of {\it multifractional} processes
such as multifractional Brownian motion
with variable amplitude (Section 6) and
multifractional $\alpha$-stable motions (Section 7).  In Section 9
we develop {\it multistable} processes, where the stability index
$\alpha(t)$ is allowed to vary.  Here the constructions are based on sums
over Poisson processes for which the required properties are
reviewed in Section 8.  In particular we
construct {\it multifractional multistable} processes, where both the local 
self-similarity index and the stability index vary.


\section{Convergence of random processes}
\setcounter{equation}{0}
\setcounter{theo}{0}
\medskip
This section is largely intended to establish notation. 
We work with two definitions of
localisability of real valued random processes, one in terms of 
 convergence of finite dimensional distributions and 
one requiring the stronger convergence in distribution, appropriate
when the sample functions are viewed as members of some metric space.

Given a probability space $(\Omega,{\cal P},\P)$, a {\it random process}
$X$ on a domain $T$ is 
a family of random variables $\{X(t): t \in T\}$.  
For our purposes $T$ will be either $\bbbr$ or a
subinterval of $\bbbr$, or sometimes  a subset of $\bbbr^{2}$
in which case we will refer to the process as a {\it random field}.

We write $X_{r}\fdd X$ to mean that a family of random processes 
$X_{r}$ converges to a process $X$ in finite-dimensional
distributions.

For processes with sample paths
in suitable function spaces, convergence in distribution
is defined in terms of a metric on the spaces.
Let $C(T)$ be the space of continuous
functions on $T \subset\bbbr$.  Writing $d^{T}(x,y) = \sup_{t\in T}|x(t)-y(t)|$
for the uniform metric on $C(T)$,
\begin{equation}
    d(x,y) = \sum_{\tau=1}^{\infty} 2^{-\tau}\min\{1,d^{[-\tau,\tau]}(x,y)\}
     \quad (x,y \in C(\bbbr))\label{met}
\end{equation}
defines a seperable metric on $C(\bbbr)$ that gives the topology 
of uniform convergence on compact subsets of $\bbbr$.

To accommodate processes with sample functions that have jumps,
let  $T$ be a closed subinterval of $\bbbr$, and let $D(T)$ denote 
the ``c\`{a}dl\`{a}g''  
functions on $T$, that is functions which are continuous on the right and
have left limits at all $t \in T$.  When $T$ is a bounded closed interval 
we define a metric $d^{[a,b]}_{S}$ on $D[a,b]$ as follows.  Let $\Phi$ be 
the class of strictly increasing continuous bijections from $[a,b]$ to
itself.  
For each $x,y \in D[a,b]$ we define
$d^{[a,b]}_{S}(x,y)$ to be the infimum of those $\delta > 0$
for which there exists $\phi \in \Phi$ such that
both $\sup_{0\leq t\leq 1}|\phi(t) - t| \leq \delta$ and 
$\sup_{0\leq t\leq 1}|x(t)-y(\phi(t) )| \leq \delta$.
Then $d^{[a,b]}_{S}$ is the
 {\it Skorohod metric} on $D([a,b])$,
see \cite[Chapter VI]{Pol1} or \cite{Bil}.  The
Skorohod metric extends to a seperable metric on $D(\bbbr)$ by
\begin{equation}
    d_{S}(x,y) = \sum_{\tau=1}^{\infty} 
    2^{-\tau}\min\{1,d^{[-\tau,\tau]}_{S}(x,y)\}
     \quad (x,y \in D(\bbbr)).\label{mets}
\end{equation}
  
Taking $T$ as $[a,b]$ or $\bbbr$, let $F(T)$ be 
either $C(T)$ or $D(T)$ with the appropriate metric as above.
Given a probability space $(\Omega,{\cal P},\P)$ 
we call $X:\Omega \to F(T)$ 
a  {\it random function} or {\it random element} of $F(T)$ 
if $X^{-1}(B) \in {\cal P}$ 
for every Borel subset $B$ of the metric space $F(T)$.  If $T'$ is
a suitable subset of $T$ and $X$ is a random function on $T$ then we may
regard the restriction of $X$ as a random function on $T'$.
When we write
$X=Y$ it will be clear from 
the context whether this refers to equality in finite dimensional
distributions or in distribution.

For $X_{r}$ and $X$ random functions in $F(T)$ where $T$ is a closed 
interval, perhaps $\bbbr$, we say that 
$X_{r}$ {\it converges in distribution} to
$X$, written $X_{r} \cd  X$, if $\E(f(X_{r})) \rightarrow \E(f(X))$
for all bounded continuous $f: F(T) \rightarrow
\bbbr$.  
Convergence in distribution is equivalent to convergence of 
finite dimensional distributions together with an appropriate stochastic 
equicontinuity condition, see \cite{Bil,Pol1}. 

Note that convergence in distribution in $C(\bbbr)$ or $D(\bbbr)$ is
equivalent to
convergence in distribution of the restrictions of the random
functions to every compact interval $[a,b]$.
A technicality here is that our functions or processes may have a domain
$U$ that is a proper interval of $\bbbr$. This presents no
difficulty, since $X_{r}$ will
generally be a
sequence of enlargements of a process about some $u$ interior to $U$, 
and the domain of definition will 
eventually include every 
bounded interval 
$[a,b]$.



\section{Localisable processes}
\setcounter{equation}{0}
\setcounter{theo}{0}
\medskip
For convenience we give the definitions of localisability at $u$ for
random processes with domain $\bbbr$, but the definitions will also
apply
in the obvious way where the domain is a real interval with $u$ as an 
interior point.  
Intuitively, a random process $Y$ on $\bbbr$ 
is localisable at $u \in \bbbr$ if it
has a unique non-trivial scaling limit at $u$. 
More precisely, we say that $Y=\{Y(t):t\in \bbbr\}$  is $h$-{\it localisable} at 
$u$ with {\it local form} the random process $Y_{u}'=\{Y_{u}'(t):t\in \bbbr\}$,
if
\begin{equation}
\frac{Y(u+rt) - Y(u)}{r^{h}} \to Y_{u}'(t)  \label{loc}
\end{equation}
as $r \searrow 0$, where 
convergence is of finite dimensional
distributions.  If $Y$ and  $Y_{u}'$ have versions in $C(\bbbr)$ or $D(\bbbr)$ and 
convergence in (\ref{loc}) is in distribution, we
say that $Y$ is {\it strongly localisable} at $u$ with {\it strong
local form} $Y_{u}'$. Of course,
strongly localisable processes are localisable, with the strong local 
form a version of the local form in $C(\bbbr)$ or $D(\bbbr)$. 
Note that the term
{\it locally asymptotically self-similar}  is sometimes used for
strong localisability.

A number of well-known processes are
$h$-localisable, in particular processes that are 
$h$-{\it self-similar}, that is $Y(rt) = r^{h}Y(t)$ for $r>0$, 
and which have {\it stationary increments} (that is $Y(t+u)$ and $Y(t)$ equal in law for $u \in 
\bbbr$).

\begin{prop}
Let $\{Y(t): t\in \bbbr\}$ be a process that is 
$h$-self-similar with stationary increments ($h$-sssi). Then $Y$ is 
 $h$-localisable at all $u\in \bbbr$ with 
$Y_{u}'= Y$.   If in addition $Y$ is in $C(\bbbr)$ or $D(\bbbr)$
then $Y$ is strongly
 $h$-localisable at all $u\in \bbbr$.
\end{prop}

\noindent{\it Proof.} If $Y$  is 
$h$-self-similar with stationary increments, then
$$\frac{Y(u+rt) - Y(u)}{r^{h}} =
\frac{Y(rt) - Y(0)}{r^{h}}=
\frac{Y(rt)}{r^{h}} = Y(t)$$
in law for all $r \neq 0$, so $Y$ is localisible at $u$.  

Further, if $Y$ is in 
$C(\bbbr)$ or $D(\bbbr)$ then 
$(Y(u+rt) - Y(u))/r^{h}$ and $Y(t)$ have identical probability
distributions, 
since probability distributions on $C(\bbbr)$ and $D(\bbbr)$ are
completely determined by their finite dimensional distributions, see
\cite{Bil}.  Thus $Y$ is strongly localisible.
\Box
\medskip

There are several important processes which are sssi
so which are strongly localisable by Proposition 3.1.

For $0<h<1$, index-$h$ fractional Brownian motion (fBm) on 
$\bbbr$ may be defined as a stochastic integral with respect to Wiener 
measure $W$:
\begin{equation}
B_{h}(t) = c(h)^{-1} \int_{-\infty}^{\infty}\left((t-x)_{+}^{h-1/2} - 
(-x)_{+}^{h-1/2}\right) W(dx), \label{fBm}
\end{equation}
where $(a)_{+}= \max\{{0,a}\}$ and $c(h)$ is a normalising constant that 
ensures that the variance $\var B_{h}(1) = 1$.  (Here, and 
throughout, we make the convention that expressions involving the 
difference of two positive parts represent an indicator function when 
the exponent is $0$, so for example, if $h=1/2$ then $(t-x)_{+}^{h-1/2} - 
(-x)_{+}^{h-1/2}$ is taken to mean $\1_{[0,t)}(x)$.) 
It is well-known \cite{EM,Fa,MV,ST} that index-$h$ fBm
 is an $h$-self-similar process  with a version in $C(\bbbr)$
that has stationary increments, so  is  
strongly localisable at all $u\in \bbbr$ with 
$(B_{h})_{u}'= B_{h}$.

The $\alpha$-stable processes form another important class of fractal 
processes of $C(\bbbr)$, or of $D(\bbbr)$ in the case of `jump'
processes, see Section 7.  Under certain 
conditions   $\alpha$-stable processes may be
sssi, see \cite[Corollary 7.3.4]{ST}, in which case by Proposition 3.1
they are strongly $h$-localisable.  

A particular instance is linear stable fractional motion:
\begin{equation}
L_{\alpha,h}(t)= \int_{-\infty}^{\infty}\left[a\left((t-x)_+^{h-1/\alpha}
-(-x)_+^{h-1/\alpha}\right) + b\left((t-x)_-^{h-1/\alpha}
-(-x)_-^{h-1/\alpha}\right) \right]M(dx),\label{lsfm}
\end{equation}
where $0<\alpha<2$ and $M$ is an $\alpha$-stable random measure
with constant skewness $\beta$ and control
measure Lebesgue measure, $0<h<1$ and $a$ and $b$ are constants,
see \cite[Section 7.4 and Chapter 10]{ST}.  
The process 
is $h$-sssi and so is $h$-localisable at all $u\in \bbbr$ with
$(L_{\alpha,h})_{u}'= L_{\alpha,h}$.  Provided that $h>1/\alpha$ 
it has a version in
$C(\bbbr)$, so is  
strongly localisable.
However, if $h<1/\alpha$ then almost surely $Y$ is unbounded on 
every interval and so is not a process of $D(\bbbr)$, though it is
nevertheless localisable.  (Note that later we will represent such processes as Poisson sums
rather than integrals with respect to random measures.)

An $\alpha$-stable L\'{e}vy motion, $0<\alpha<2$ is a process in $D(\bbbr)$
with stationary independent increments which have a strictly $\alpha$-stable
distribution.   It may be represented as
\begin{equation}
L_{\alpha}(t)= M([0,t]) \label{aslm}
\end{equation}
where $M$ is an $\alpha$-stable random measure on $\bbbr$ with constant
skewness intensity, see \cite[Section 7.5]{ST}.  
Then $L_{\alpha}$ is $1/\alpha$-sssi, 
and so is strongly $1/\alpha$-localisable.

In later sections we will give general constructions of localisable processes
where the local form $Y_{u}'$ varies with $u$.
For now we note that localisability behaves well under reasonably smooth
changes of coordinates.  In particular 
the following proposition allows the introduction of varying `local
amplitude' for localisable processes.

\begin{prop}
Let  $U$ be an interval with $u$ an interior point.
Suppose that $\{Y(t): t\in U\}$ is $h$-localisable 
(resp. strongly $h$-localisable) at $u$.
Let 
$a:U\to\bbbr$ satisfy an $\eta$-H\"{o}lder condition on $U$, that is
$$ |a(t)-a(t')| \leq c |t-t'|^{\eta} \quad (t,t' \in U),$$
where $\eta>h$.   Then $aY= \{ a(t)Y(t): t \in U\}$
is $h$-localisable (resp. strongly $h$-localisable) at $u$ with
$(aY)_{u}' = a(u)Y_{u}'$.
\end{prop}

\medskip
\noindent{\it Proof.} We have
$$\frac{a(u+rt)Y(u+rt) - a(u)Y(u)}{r^{h}} =
a(u+rt) \frac{Y(u+rt)- Y(u)}{r^{h}}+
Y(u)\frac{a(u+rt) - a(u)}{r^{h}}.$$
The result now follows on letting $r \to 0$ with 
the appropriate form of convergence, noting that the right-hand term
has zero limit almost surely.
\Box
\medskip


\section{Localisable processes with prescribed local form}
\setcounter{equation}{0}
\setcounter{theo}{0}
\medskip
We aim to construct localisable funtions with prescribed local form
by `joining together' localisable processes $\{X(t,v): t \in
U\}$ over a range of $v$.  Thus we 
 seek conditions that ensure $Y=\{X(t,t) : t\in U\}$
looks locally like $\{X(t,u): t \in U\}$ when $t$ is
close to $u$.

Let $U$ be an interval with $u$ an interior point.  Let 
$\{X(t,v) : (t,v) \in U\times  U\}$ 
be a random field and  
let $Y$ be the diagonal process
$Y=\{X(t,t) : t\in U\}$. 
We want $Y$ and $X(\cdot,u)$ to have the same local
forms at $u$, that is
$Y_{u}'(\cdot) = X_{u}'(\cdot,u)$ where $X_{u}'(\cdot,u)$ is the local
form of $X(\cdot,u)$ at $u$.  Thus we require
\begin{equation}
\frac{X(u+rt,u+rt) - X(u,u)}{r^{h}} \fdd X_{u}'(t,u) \label{locform}
\end{equation}
as $r \searrow 0$.  The following theorem gives a sufficient
condition for this to occur.

\medskip
\begin{theo}
Let $U$ be an interval with $u$ an interior point. 
Suppose that for some $0<h<\eta$ the process  $\{X(t,u) : t \in 
U\}$  
is $h$-localisable at $u\in U$
with local form $X_{u}'(\cdot,u)$ and
\begin{equation}
\P(|X(v,v) - X(v,u)|\geq |v-u|^{\e}) \to 0 \label{cond}
\end{equation}
as $v\to u$.  Then $Y=\{X(t,t) : t\in U\}$ is $h$-localisable at $u$
with $Y_{u}'(\cdot) = X_{u}'(\cdot,u)$.

In particular, this conclusion holds if for some $p>0$ and $\eta>h$
\begin{equation}
\E(|X(v,v) - X(v,u)|^{p})=O(|v-u|^{\e p}) \label{cond1}
\end{equation}
as $v\to u$.
\end{theo}

\medskip
\noindent{\it Proof.}
For $r\neq 0$
\begin{eqnarray}
    \frac{Y(u+rt) - Y(u)}{r^{h}}
    & = &\frac{X(u+rt,u+rt) - X(u,u)}{r^{h}}\nonumber  \\
     & & \hspace{-3cm}= \quad \frac{X(u+rt,u+rt) - X(u+rt,u)}{r^{h}}
     + \frac{X(u+rt,u) - X(u,u)}{r^{h}}.\label{B}
\end{eqnarray}
Fix $t \in \bbbr$ and $c >0$.  Let $r_{0}$ be sufficiently
small to ensure that if $0<r<r_{0}$ then both $u \pm rt \in U$ and
$c r^{h} \geq (r|t|)^{\e}$.  Then for $0<r<r_{0}$
\begin{eqnarray*}
     & & \hspace{-2.5cm}
     \P\left( \frac{|X(u+rt,u+rt) - X(u+rt,u)|}{r^{h}} \geq c
     \right)  \\
     & \leq & \P\left(|X(u+rt,u+rt) - X(u+rt,u)|\geq (r|t|)^{\e}\right)  \\
     & \leq & \P \left(|X(u+rt,u+rt) - X(u+rt,u)|\geq  |(u+rt)-u|^{\e} \right)
     \to 0
\end{eqnarray*}
as $r \searrow 0$, by (\ref{cond}).  Thus for all $t \in \bbbr$,
$$\frac{X(u+rt,u+rt) - X(u+rt,u)}{r^{h}} \to 0$$
in probability  and so in finite
dimensional distributions.  Moreover,
$$\frac{X(u+rt,u) - X(u,u)}{r^{h}} \fdd X_{u}'(t,u),$$
since $X(\cdot,u)$ is localisable at $u$. We conclude from (\ref{B})
that $Y$ is localisable at $u$ with local form
$Y_{u}'(\cdot) = X_{u}'(\cdot,u)$.

If (\ref{cond1}) holds,
Markov's inequality implies (\ref{cond}) (with $\eta$ replaced by 
some $h<\eta'<\eta$) and the
conclusion follows.
\Box

\medskip

Although Theorem 4.1 is valid for all $h>0$, it is normally
applied
 with $0<h<1$.  If $X(\cdot,u)$ is $h$-localisable for $h>1$ then the 
limit of (\ref{B}) is usually dominated by the left-hand term giving
that $Y$ is $1$-localisable, see Theorem 9.4 for an example of this. 

\section{Strongly localisable processes with prescribed local form}
\setcounter{equation}{0}
\setcounter{theo}{0}
\medskip
We obtain an analogue of Theorem 4.1 in the strongly
localisable case, that is a criterion for convergence  in
distribution in (\ref{loc}).

\medskip
\begin{theo}\label{slcp}
Let $F(\bbbr)$ be either $C(\bbbr)$   endowed
with the metric $d$ or $D(\bbbr)$ with $d_{S}$, see $(\ref{met})$ or $(\ref{mets})$.
Let $U$ be an interval with $u$ an interior point.  Suppose that
for some $h>0$ the process $\{X(t,u): t \in U\}$ of $F(U)$
is strongly $h$-localisable at $u$, with 
local form $X_{u}'(\cdot,u)$ a random function of $F(\bbbr)$.   
Suppose that
 for all
$c> 0$ 
\begin{equation}
\P\left(\sup_{0<|v-u|<\epsilon}
\frac{|X(v,v) - X(v,u)|}{|v-u|^{h}}> c\right)
\to 0 \label{cond3}
\end{equation}
as $\epsilon \to 0$. 
If the process $Y=\{X(t,t) : t\in U\}$ is in $F(U)$ then $Y$ 
is strongly $h$-localisable at $u$ 
 with $Y_{u}'(\cdot) = X_{u}'(\cdot,u)$.
 
 In particular,  this conclusion holds if for some
 $\eta >h $ we have
\begin{equation}
\sup_{v \in U, v\neq u}
\frac{|X(v,v) - X(v,u)|}{|v-u|^{\e}} 
<\infty \quad \label{conda}
\end{equation}
almost surely. 
\end{theo}

\medskip
\noindent{\it Proof.}
First consider $(C(\bbbr),d)$.   For each positive $\tau$ and $r$
sufficiently small,
\begin{eqnarray}
    \P\Bigg(\sup_{0<|t| \leq \tau}& & \hspace{-1cm}
    \frac{|X(u+rt,u+rt) - X(u+rt,u)|}{r^{h}}> c\Bigg)\nonumber  \\
    & \leq & \P\left(\sup_{0<|t| \leq \tau}
    \frac{|X(u+rt,u+rt) - X(u+rt,u)|}{|rt|^{h}}> \frac{c}{\tau^{h}}\right)\nonumber  \\
     &\leq &\P\left(\sup_{0<|v-u| \leq r\tau}
    \frac{|X(v,v) - X(v,u)|}{|v-u|^{h}}> \frac{c}{\tau^{h}}\right) \to 0
     \label{B1}
\end{eqnarray}
as $r\to 0$. Thus, the restriction of 
$\displaystyle{\frac{X(u+rt,u+rt) - X(u+rt,u)}{r^{h}}}$ to $[-\tau,\tau]$, 
 converges to $0$ in
probability in $(C[-\tau,\tau],d^{[-\tau,\tau]})$ as $r \to 0$.  From the definition 
(\ref{met}) of $d$,
convergence in probability on every 
bounded interval implies 
convergence in probability on $(C(\bbbr),d)$, so
\begin{equation}
    \frac{X(u+rt,u+rt) - X(u+rt,u)}{r^{h}}\cp 0 \label{cp}
\end{equation}
    in $(C(\bbbr),d)$.
Then
\begin{eqnarray}
    \frac{Y(u+rt) - Y(u)}{r^{h}}
    & = &\frac{X(u+rt,u+rt) - X(u,u)}{r^{h}}\nonumber  \\
     &= & \frac{X(u+rt,u+rt) - X(u+rt,u)}{r^{h}}
     + \frac{X(u+rt,u) - X(u,u)}{r^{h}} \nonumber  \\
     & \cd & X_{u}'(t,u),
     \label{B2}
\end{eqnarray}
as $r\to 0$,  since $X$ is localisable at $u$.  Here we use 
a standard property \cite[Theorem 4.1]{Bil}, that for random elements 
$Z_{r},Z,W_{r},W$ of some metric space $(M,\rho)$,  if $Z_{r}\cd Z$ and 
$\rho(Z_{r},W_{r}) \cp 0$, then $W_{r}\cd Z$.  

Turning to $(D(\bbbr),d_{S})$, if $X(t,u) \in D(\bbbr)$ and 
(\ref{cond3}) holds, the same argument using (\ref{B1}) implies 
convergence in probability in (\ref{cp}) with respect to the metric 
$d_{S}$. (Note that $d_{S}^{[-\tau,\tau]}(f,0) \leq \sup_{t \in [-\tau,\tau]} |f(t)|$ for 
each $\tau$ for $f\in D(\bbbr)$.)  Convergence in distribution then 
follows just as for $C(\bbbr)$.

Finally, (\ref{cond3}) is an immediate consequence of (\ref{conda}) if $h< 
\eta$.
\Box
\medskip

To utilise Theorem 5.1 we need to verify (\ref{conda}), that is to show
that
$Z(v) = (X(v,v)-X(v,u))/|v-u|^{\eta}$ 
is bounded as $v$ ranges across an interval.  The following form of 
Kolmogorov's continuity theorem will be extremely useful for this.

\medskip
\begin{theo} (Kolmogorov's continuity theorem)
 Let $\{Z(v) : v \in T\}$ be a random process where $T$ is a bounded
 subset of  $\bbbr^{n}$.  If for some $p>0,\epsilon>0$ and $c>0$
 $$\E|Z(v) - Z(v')|^{p} \leq c |v-v'|^{n+\epsilon}\quad (v,v' \in
 T),$$
 then $Z$ has a continuous version that is almost surely
 $\eta$-H\"{o}lder continuous for all
 $0 <\eta <\epsilon/p$.
 \end{theo}
 
\medskip
\noindent{\it Proof.}  See, for example,  \cite[Theorem 25.2]{RW}.
\Box


\section{\!Multifractional Brownian motion with variable amplitude}
\setcounter{equation}{0}
\setcounter{theo}{0}
\medskip
A number of constructions of multifractional Brownian motion, 
a process with index-$h(u)$ fractional Brownian
motion as its local form at $u$, have been
given, see \cite{AA,AL2,BJR,PL}.   To demonstrate our method we 
indicate briefly a
straightforward construction of multifractional Brownian motion, that
is strongly localisable 
with a given local index and amplitude.   

As in \cite{PL} we model our definition on (\ref{fBm}) but allow $h$ to vary.  By
virtue of  Proposition 3.2 variable local amplitude
presents no difficulty.  Let $U$ be a bounded closed interval and let
$h:U \to (0,1)$ satisfy an $\eta$-H\"{o}lder condition
\begin{equation}
|h(v)-h(v')| \leq k|v-v'|^{\eta} \quad (v,v' \in U) \label{hol}
\end{equation}
where $0<\eta \leq 1$.

\medskip
\begin{theo} (Multifractional Brownian motion)
Let $u \in \bbbr$ and let $U$ be  a closed interval with $u$ an
interior point.  Suppose
that  $h: U \to (0,1)$ and $a: U \to \bbbr^{+}$ 
both satisfy an $\eta$-H\"{o}lder condition where $h(u) < \eta \leq 1$.
Define
\begin{equation}
Y(t) = a(t)\int_{-\infty}^{\infty}\left((t-x)_{+}^{h(t)-1/2} - 
(-x)_{+}^{h(t)-1/2}\right) W(dx) \quad (t \in U). \label{mbmc}
\end{equation}
Then $Y$ is strongly $h(u)$-localisable at $u$ with
$Y_{u}' = a(u)c(h(u))B_{h(u)}$ where
$B_{h}$ is index-$h$ fBm and where $c(h)$ is the normalisation constant
in $(\ref{fBm})$.
\end{theo}

\noindent{\it Proof.}  
By Proposition 3.2 it is enough to consider the case where $a(v)\equiv
1$.
We define a random field by the
stochastic integral
\begin{equation}
X(t,v) = \int_{-\infty}^{\infty}\left((t-x)_{+}^{h(v)-1/2} - 
(-x)_{+}^{h(v)-1/2}\right) W(dx) \quad (t,v \in U), \label{mBm}
\end{equation}
where $W$ is Wiener measure on $\bbbr$. Since the integrand of (\ref{mBm}) is square
integrable, $X(t,v)$ exists a.s. with mean $0$ for all $t,v \in U$. A 
mean value estimate applied to (\ref{mBm}) easily gives that 
$\E (X(t,v)-X(t',v'))^{p} \leq  c(|t-t'|^{p\eta} +|v-v'|^{p\eta})$,  
first for $p=2$ and then for all $p>0$ since the increments are 
Gaussian. By Kolmogorov's criterion, for all $\epsilon >0$ such that $h(u)
<\eta-\epsilon$ there is an a.s. finite random variable $C$ such
that
$$|X(v,v)-X(v,u)| \leq C |v-u|^{\eta-\epsilon} \quad (v \in U)$$
so  (\ref{conda}) holds (with $\eta$ replaced by $\eta - \epsilon$).  
But $X(\cdot,u)= B_{h(u)}(\cdot)$ which
is sssi so is strongly $h(u)$-localisable at $u$ by
Proposition 3.1. 
Theorem 5.1 implies that $Y=\{X(t,t): t\in T\}$ 
is strongly  $h(u)$-localisable at 
$u$ with $Y_{u}'(\cdot) = X_{u}'(\cdot,u)=(B_{h(u)})_{u}'(\cdot)=B_{h(u)}(\cdot)$.
\Box
\medskip


\section{Multifractional stable processes}
\setcounter{equation}{0}
\setcounter{theo}{0}
\medskip

Multifractional Brownian motion generalizes fractional Brownian motion 
by allowing the parameter $h$ to  vary with time. By
working with a stochastic integral with 
respect to an $\alpha$-stable measure instead of Wiener measure, 
we now construct multifractional stable processes  with the
local scaling exponent depending on $t$.

Recall that a process $\{X(t):t\in T\}$, where $T$ is
generally a subinterval of $\bbbr$,
is called $\alpha$-{\it stable} $(0<\alpha \leq 2)$ 
if all its finite-dimensional distributions are
$\alpha$-stable, see the encyclopaedic work on stable 
processes \cite{ST}. Note that $2$-stable processes are just Gaussian 
processes. 

Many stable processes admit a stochastic integral representation.
Write $S_{\alpha}(\sigma,\beta,\mu)$ for the $\alpha$-stable distribution
with scale parameter $\sigma$, skewness $\beta$ and shift-parameter 
$\mu$; we will assume throughout that $\mu=0$.  Let $(E,{\cal E},m)$ be a
sigma-finite measure space  (which will be Lebesgue measure in our
examples).  
Taking $m$  as the control measure and $\beta: E \to [-1,1]$ a
measurable function, this defines an
$\alpha$-stable random measure $M$ on $E$ such that
for $A\in {\cal E}$ we have that 
$M(A) \sim S_{\alpha}(m(A)^{1/\alpha},
\int_{A}\beta(x)m(dx)/m(A),0)$.  If $\beta =0$ then the process is
{\it symmetric} $\alpha$-stable or S$\alpha$S.

Let
$$    {\cal F}_{\alpha}\equiv {\cal F}_{\alpha}(E,{\cal E}, m) 
= \{ f: f \mbox{ is measurable and } \|f\|_{\alpha} < \infty\},$$
where $\|\,\|_{\alpha}$ is the quasinorm (or norm if $1<\alpha \leq
2$) given by 
\begin{equation}
\|f\|_{\alpha}     
=\left\{ 
\begin{array}{cc}
    \left( \int_E |f(x)|^{\alpha}m(dx)\right)^{1/\alpha} & (\alpha\neq 1)  \\
     \int_E |f(x)|m(dx) 
     + \int_E |f(x)\beta(x) \ln |f(x)| | m(dx)  & (\alpha = 1)
\end{array}
 \right. 
    \label{normdef}
\end{equation}
The stochastic
integral of $f\in  {\cal F}_{\alpha}(E,{\cal E}, m)$ with respect to
$M$ then exists  \cite[Chapter 3]{ST} with
\begin{equation}
I(f) =\int_E f(x)M(dx)\sim
S_{\alpha}(\sigma_{f},\beta_{f},0),\label{alint}
\end{equation}
where
$$\sigma_{f}=\|f\|_{\alpha}, \quad
\beta_{f} = \frac{\int f(x)^{<\alpha>}\beta(x) m(dx)}
{\|f\|_{\alpha}^{\alpha}},$$
writing $a^{<b>} \equiv \mbox{sign}({a}) |a|^{b}$, see \cite[Section
3.4]{ST}.
In particular,
\begin{equation}
\E|I(f)|^{p} = \left\{
\begin{array}{cc}
   c(\alpha,\beta,p)\|f\|_{\alpha}^{p}  &  (0<p<\alpha)  \\
   \infty  & (p \geq \alpha)
\end{array}\right. \label{alintexp}
\end{equation}
where $c(\alpha,\beta,p)<\infty$, see \cite[Property 1.2.17]{ST}.

When $0<\alpha<1$ there is a non-negative stable subordinator measure $M'$
associated with $M$ so that $M'(A) \sim S_{\alpha}(m(A)^{1/\alpha},
1,0)$.  In particular, for $f \in {\cal F}_{\alpha}$, 
\begin{equation}
|I(f)| \leq \int_E |f(x)|M' (dx).  \label{ss}
\end{equation}

We will be concerned with processes that may be
expressed as stochastic integrals
\begin{equation}
X(t)=\int_E f(t,x)M(dx) + \mu(t), 
\quad (t \in T),\label{intrep}
\end{equation}
where $f(t,\cdot)$ is a jointly measurable family of functions in
${\cal F}_{\alpha}(E,{\cal E}, m)$
and $\mu(t)$ are real numbers.  Note
that if $\mbox{\rm esssup}_{a \leq t \leq b} f(t,x)= \infty$ for all $x \in A$
for some
$A \subset E$ with $m(A)>0$ then $X(t)$ will be unbounded a.s. on
the interval $[a,b]$, see \cite[Section 10]{ST}. 

Here we consider the localisibility at $u$ of processes 
defined in terms of random fields
\begin{equation}
X(t,v)= \int_E f(t,v,x)M(dx) + \mu(t,v) \quad(t,v \in U) \label{field}
\end{equation}
where $f(t,v,.)\in {\cal F}_{\alpha}$ and $\mu(t,v) \in \bbbr$ for all
$t, v \in U$ for some interval $U$. We assume throughout that  $f(t,v,x)$ is
measurable on $U\times U \times E$.

The term $\mu(t,v)$ is  easily dealt with: if  $v \mapsto \mu(v,v)$ 
is pointwise $\eta$-H\"older at $v = u$, 
that is
$|\mu(v,v) - \mu(u,u)| \leq k |u-v|^{\eta}$
for $v $ close to $u$, where $0<h < \eta \leq 1$, then the
$h$-localisability of $Y=\{X(t,t): t\in U\}$ at $u$ and its local form are
unaffected if we set $\mu(t,v) =0$, so we assume this
throughout this section.

The following proposition gives conditions for  $Y$  to have a
continuous or bounded version, which is needed for strong
localisability to be meaningful.  Note that these sufficient
conditions are geared towards our context; for other aspects see \cite[Chapters
10,12]{ST}.

\begin{prop}
Let $U$ be a closed interval.
Let $X$ be a random field defined by 
\begin{equation}
X(t,v)= \int_{E} f(t,v,x)M(dx) \quad (t,v \in U) \label{fieldcont}
\end{equation}
where $f(t,v,\cdot)\in {\cal F}_{\alpha}$ are jointly measurable
and $M$ is an $\alpha$-stable
random measure with control measure $m$ and measurable skewness.

$(a)$  Let $0<\alpha<1$.  If
\begin{equation}
\| \sup_{t,v \in U}|f(t,v,x)|\|_{\alpha}< \infty, \label{supfin}
\end{equation}
 then the random field  $(\ref{fieldcont})$ has a bounded version.
 
 If in addition $\{f(t,v,x) : x \in E\}$ is an equiuniformly
 continuous family
 for $t,v \in U$, then $(\ref{fieldcont})$ has a continuous version.
 
$(b)$ Let  $1<\alpha<2$ and  $1/\alpha<\eta\leq 1 $.  If
\begin{equation}
    \|f(t,v,\cdot)-f(t',v',\cdot)\|_{\alpha} 
    \leq k\left(|v-v'|^{\eta}+|t-t'|^{\eta}\right)
    \quad (t,t',v,v' \in U),
    \label{cty1}
\end{equation}
then $Y=\{X(t,t): t\in U\}$ has a continuous version for $t \in U$,
satisfying an a.s. $\beta$-H\"{o}lder condition for all 
$ 0 <\beta <(\eta\alpha-1)/\alpha$.
\end{prop}

 \noindent{\it Proof.} 
 (a)  Since $0<\alpha<1$ there exists a 
 stable subordinator measure $M'$
associated with $M$, so that $M'$ has control measure $m$ and
$M'(A) \sim S_{\alpha}(m(A)^{1/\alpha},
1,0)$.  By (\ref{ss}), for $t,v \in U$,
\begin{equation*}
    |X(t,v)|\leq \int |f(t,v,x)|M'(dx)
 \leq \int \sup_{t,v \in U}|f(t,v,x)|M'(dx)
\equiv Z,
\end{equation*}
where $Z$ is an almost surely finite random variable by (\ref{supfin}),
so $X(t,v)$ is a.s. bounded for $t,v \in U$.

Now assume also the equicontinuity condition.  Given $\epsilon>0$ we
may, since $E$ is $\sigma$-finite, choose $D\subset E$ such that  
$\int_{E\setminus D}\left( \sup_{t,v \in U}|f(t,v,x)|\right)^{\alpha}m(dx)<
\epsilon^{\alpha}$.  By equiuniform continuity we may find $\delta>0$ such that for
all $x \in E$ and $|(t,v)-(t',v')|<\delta$ we have $|f(t,v,x)-f(t',v',x)|<
m(D)^{-1/\alpha}\epsilon$.   Then if $|(t,v)-(t',v')|<\delta$,  (\ref{ss})
gives
\begin{eqnarray*}
    |X(t,v)-X(t',v')|
    &\leq& \int_{E} |f(t,v,x)-f(t',v',x)|M'(dx) \\
    &\leq &2\int_{E\setminus D} \sup_{t,v \in U}|f(t,v,x)|M'(dx) 
        + \int_{D} \frac{\epsilon}{m(D)^{1/\alpha}}M'(dx) \equiv
	Z_{\epsilon},
\end{eqnarray*}
say, where $Z_{\epsilon}$ is a random variable.  Fix $0<p<\alpha$.  
By (\ref{alintexp}) there is a constant 
$c$ independent of $\epsilon$ such that
$$\E|Z_{\epsilon}|^{p} \leq c \epsilon^{p}.$$
Thus choosing $\epsilon(n)\, (n=1,2,\ldots)$ such that 
$\E|Z_{\epsilon(n)}|^{p} \leq 2^{-n}$, there are corresponding  
$\delta_{n}$ such that
$$\sup_{|(t,v)-(t',v')|<\delta_{n}}|X(t,v)-X(t',v')|
\leq  Z_{\epsilon(n)}.$$
Since $\sum_{n=1}^{\infty}\E|Z_{\epsilon(n)}|^{p} <\infty$, the Borel-Cantelli
lemma gives that $Z_{\epsilon(n)} \to 0$ almost surely, so 
$\sup_{|(t,v)-(t',v')|<\delta_{n}}|X(t,v)-X(t',v')|\to 0$ a.s. as 
$n \to \infty$, giving  continuity of $X(t,v)$ a.s.
 
 (b) From (\ref{fieldcont})
$$X(t,v)-X(t',v')  =\int  
\left(f(t,v,x) - f(t',v',x) \right)M(dx).$$
This integrand is in ${\cal F}_{\alpha}$,  so for $0<p<\alpha$,
estimate (\ref{alintexp}) gives
\begin{eqnarray*}
   \E |X(t,v)-X(t',v')|^{p}  & \leq & c_{1}
     \left\|f(t,v,\cdot) - f(t',v',\cdot) \right\|_{\alpha}^{p}  \\
     & \leq & c_{2}\left(|v-v'|^{\eta p}+|t-t'|^{\eta p}\right)
\end{eqnarray*}
by (\ref{cty1}) where $c_{1}$ and $c_{2}$ are independent of $t,t',v,v' \in U$.  
Specialising, 
\begin{equation*}
   \E |Y(t)- Y(t')|^{p}=\E |X(t,t)- X(t',t')|^{p}  
   \leq 2c_{2}|t-t'|^{\eta p} 
\end{equation*} 
for $t,t' \in U$.  

Since $\eta > 1/\alpha$ we may choose $0<p<\alpha$ 
such that $\eta p >1$.   Kolmogorov's Theorem 5.2 gives that 
$Y$ has a continuous version  for $t \in U$ with an a.s.
$\beta$-H\"{o}lder condition for all 
$ 0 <\beta <(\eta p-1)/p$ for all $p<\alpha$. 
 \Box
\medskip

We require the following calculus lemma.

\begin{lem}
Let $U$ be an interval and let $f:U \to \bbbr$ be  
continuously differentiable with $f'$ satisfying an $\eta$-H\"{o}lder
condition
\begin{equation}
|f'(v)-f'(w)| \leq k|v-w|^{\eta} \quad (v,w \in U) 
\label{diffhol}
\end{equation}
for some $0<\eta \leq 1$. Let $v,w,u \in U$ with $ v\neq u, w\neq u$.  Then 
\begin{equation}
\left| \frac{f(v)-f(u)}{v-u} - \frac{f(w)-f(u)}{w-u}\right|
 \leq 2^{\eta} k |v-w|^{\eta} .\label{cal}
\end{equation}
\end{lem}

\noindent{\it Proof.}  
We may assume without loss of generality 
that $v<w$ and $u<w$. Write $\displaystyle{g(v) =  \frac{f(v)-f(u)}{v-u}}$.
We consider three cases.

(a) If $v<u<w$, then by the mean value theorem there exist $v_{0} \in
(v,u)$ and $w_{0} \in (u,w)$ such that $g(v) = f'(v_{0})$ and 
$g(w) = f'(w_{0})$.  Then
$$ |g(v)-g(w)| = |f'(v_{0})-f'(w_{0})|\leq k |v_{0}-w_{0}|^{\eta} 
\leq k |v-w|^{\eta} .$$

(b) If $u<v<w$ and  $|w-v| \geq |v-u|$, then  
$|w-v| \geq \frac{1}{2}|w-u|$.   There exist $v_{0} \in
(u,v)$ and $w_{0} \in (u,w)$ such that $g(v) = f'(v_{0})$ and 
$g(w) = f'(w_{0})$, so 
$$ |g(v)-g(w)| = |f'(v_{0})-f'(w_{0})|\leq k |v_{0}-w_{0}|^{\eta} 
\leq k |w-u|^{\eta} \leq k 2^{\eta}|w-v|^{\eta} .$$

(c) If $u<v<w$ and  $|w-v| \leq |v-u|$, we apply the mean value theorem
to $g$.  Thus there exists $s \in (v,w)$ such that
\begin{eqnarray*}
  g(v) -g(w)  & = & (v-w) g'(s)  \\
     & = &  (v-w)\frac{(s-u)f'(s) - f(s)+f(u)}{(s-u)^{2}}  \\
     & = & (v-w)\frac{f'(s) - f'(z)}{(s-u)} 
\end{eqnarray*}
where $z \in (u,s)$ using the mean value theorem again. Hence
\begin{eqnarray*}
  |g(v) -g(w)|  & \leq & k\frac{|v-w||s-z|^{\eta}}{|s-u|}  \\
     & \leq & k|v-w||s-u|^{\eta-1}  \\
     & \leq & k|v-w||v-u|^{\eta-1} \\
     & \leq & k|v-w|^{\eta}.
\end{eqnarray*}
\Box
\medskip

The following theorem gives conditions that allow the transfer of localisability 
properties from $X(\cdot,u)$ to $Y= \{ X(t,t): t\in U\}$ in the $\alpha$-stable
case, generalising the results of Section 6 in the Gaussian case.

\begin{theo}
Let $U$ be a closed interval with $u$ an interior point.  
Let $X$ be a random field defined by 
\begin{equation}
X(t,v)= \int f(t,v,x)M(dx) \quad (t,v \in U) \label{field1}
\end{equation}
where $f(t,v,\cdot)\in {\cal F}_{\alpha}$ are jointly measurable
and $M$ is an $\alpha$-stable
random measure with control measure $m$ and measurable skewness.

$(a)$ Suppose that $0<\alpha \leq 2$ and the process $X(\cdot,u)$ 
is $h$-localisable at $u$ with  $h>0$. 
    Suppose that for some $\eta >h$
\begin{equation}
    \|f(t,v,\cdot) - f(t,u,\cdot)\|_{\alpha} \leq k_{1} |v-u|^{\eta}
    \quad (t,v \in U). \label{hol2}
    \end{equation}
Then $Y =\{ X(t,t): t\in U\}$ is $h$-localisable at $u$ with
local form $Y'_u(\cdot)=X'_u(\cdot,u).$

$(b)$ Suppose that  $0<\alpha<1$ and that  $X(\cdot,u)$ is strongly 
$h$-localisable in $C(\bbbr)$ (resp. $D(\bbbr)$)
at $u$.  Suppose that for some $\eta>h$ 
\begin{equation}
    |f(t,v,x) - f(t,u,x)| \leq k_{1}(x) |v-u|^{\eta}
    \quad (t,v \in U, x \in E),\label{hol4}
\end{equation}
 where $k_{1}(\cdot)\in{\cal F}_{\alpha}$.
 If $Y =\{ X(t,t): t\in U\}$ has a version in $C(U)$ (resp. $D(U)$) (see
 Proposition 7.1$($a$)$),
 then $Y$ is strongly $h$-localisable at $u$ 
 in $C(\bbbr)$ (resp. $D(\bbbr)$) with
$Y'_u(\cdot)=X'_u(\cdot,u).$

$(c)$ Suppose that  $1<\alpha \leq 2$, that  $\eta>1/\alpha $ and that
$X(\cdot,u)$ 
is strongly $h$-localisable in $C(\bbbr)$ or $D(\bbbr)$
at $u$.  Suppose that for all $t,v \in U$ the partial derivative
$f_{v}(t,v,\cdot) \in {\cal F}_{\alpha}$ with
\begin{equation}
    \left|f_{v}(t,v,x)
    - f_{v}(t,v',x)\right| 
    \leq k_{1}(t,x) |v-v'|^{\eta}
    \quad (t,v,v' \in U, x \in E),
    \label{hol5}
\end{equation}
where $\sup_{t \in U}\|k_{1}(t,\cdot)\|_{\alpha}< \infty$, and that
\begin{equation}
    \sup_{v \in U}\left|f_{v}(t,v,x)
    - f_{v}(t',v,x)\right|
    \leq k_{2}(t,t',x)
    \quad (t,t' \in U, x \in E),
    \label{hol7}
\end{equation}
where $\|k_{2}(t,t',\cdot)\|_{\alpha}\leq c|t-t'|^{\eta}$.
Then $Y =\{ X(t,t): t\in U\}$ is strongly $h$-localisable at $u$ 
in $C(\bbbr)$ with
$Y'_u(\cdot)=X'_u(\cdot,u)$.
\end{theo}

\medskip
\noindent{\it Proof.}
(a) We have
\begin{equation}
    X(t,v)-X(t,u) = \int \left(f(t,v,x)- f(t,u,x) \right)M(dx) 
 \label{diff}
\end{equation}
so, taking $0<p<\alpha$ and using (\ref{alintexp}), there is a
constant $c_{1}$ such that
\begin{eqnarray*}
\E |X(t,v)-X(t,u)|^{p} &\leq & c_{1} \left\|f(t,v,\cdot)-
f(t,u,\cdot)\right\|_{\alpha}^{p}\\
&\leq & c_{1}k_{1}|v-u|^{\eta p}. 
\end{eqnarray*}
The conclusion  follows from
Theorem 4.1. 

(b) Since $0<\alpha<1$ there exists a 
 stable subordinator measure $M'$
associated with $M$, so that $M'$ has control measure $m$ and
$M'(A) \sim S_{\alpha}(m(A)^{1/\alpha},
1,0)$.  Applying (\ref{ss}) to (\ref{diff}) and using  (\ref{hol4}), 
gives that for $t,v \in U$
\begin{eqnarray*}
    |X(t,v)-X(t,u)|&\leq &\int |f(t,v,x) - f(t,u,x)|M'(dx) \\
&\leq &|v-u|^{\eta}\int k_{1}(x)M'(dx) \\
&\leq & |v-u|^{\eta} Z,
\end{eqnarray*}
where $Z$ is an a.s. finite random variable.
Thus (\ref{conda}) holds and Theorem 5.1 gives that
 $Y$ is strongly localisable at $u$.

(c)  It is easy to check that $Y$ satisfies the conditions of Theorem 7.1(b)
and so has a continuous version.
Again we verify (\ref{conda}).
Define
\begin{eqnarray*}
    Z(t,v) &=&\frac{X(t,v)-X(t,u)}{v-u}\quad (t,v \in U, v\neq u)\\
    &=& \int g(t,v,x) M(dx)
\end{eqnarray*}
where 
$$g(t,v,x) = \frac{f(t,v,x) - f(t,u,x)}{v-u}. $$
Applying  Lemma 7.2 with $f(v)= f(t,v,x)$ and noting (\ref{hol5}), we 
get
\begin{equation}
    |g(t,v,x)- g(t,v',x)|\leq 2^{\eta}k_{1}(t,x)
    |v-v'|^{\eta}.\label{gtvx}
\end{equation}
Also
\begin{eqnarray}
    |g(t,v,x)- g(t',v,x)|
    &= &\frac{1}{|v-u|}\left|(f(t,v,x)  - f(t',v,x)) -  (f(t,u,x) - 
    f(t',u,x))\right| \nonumber\\
&\leq &\left|f_{v}(t,v_{1},x)
-  f_{v}(t',v_{1},x)  
\right|\nonumber\\
&\leq & k_{2}(t,t',x), \label{gtvx2}
\end{eqnarray}
for some $v_{1} \in (u,v)$, on applying the mean value theorem to $f(t,v,x)  - 
f(t',v,x)$.  From (\ref{gtvx}) and (\ref{gtvx2}) together
with the conditions on $k_{1}$ and $k_{2}$ we get
$$\|g(t,v,\cdot)- g(t',v',\cdot)\|_{\alpha} \leq c_{1}(|v-v'|^{\eta}+
|t-t'|^{\eta}).$$

Applying Proposition 7.1(b) to $\{Z(v,v): v \in U\}$, it follows that 
$ Z(v,v) =(X(v,v) - X(v,u))/(v-u)$ has a version that is 
a.s. continuous and bounded for $v \in U$.  Thus (\ref{conda}) holds and strong 
localisability follows from Theorem 5.1.
 \Box
\medskip

We illustrate Theorem 7.3 by constructing 
processes whose local forms are 
linear stable fractional motions $L_{\alpha,h(t)}$, 
see (\ref{lsfm}). Overlapping results with a different emphasis are
given in \cite{ST1,ST2}.
The following process is termed a {\it linear stable multifractional motion}:
\begin{align}
Y(t)= \int_{-\infty}^{\infty}\Big[a & \left((t-x)_+^{h(t)-1/\alpha}
-(-x)_+^{h(t)-1/\alpha}\right) \nonumber \\
& + b\left((t-x)_-^{h(t)-1/\alpha}
-(-x)_-^{h(t)-1/\alpha}\right) \Big]M(dx)
\quad (t \in \bbbr),\label{lsmmdef}
\end{align}
where $M$ is an $\alpha$-stable random measure ($0< \alpha <2$)
with constant 
skewness intensity $\beta$ and control measure Lebesgue measure, 
with $h(t) \in (0,1)$ for all $t \in \bbbr$, and $a$ and $b$ real 
numbers. (Recall that $(w)_+=\max\{0,w\}$ and 
$(w)_- = -(w)_+$ for $w \in 
\bbbr$.)

To investigate localisability, we introduce the random field 
\begin{align}
X(t,v)= \int_{-\infty}^{\infty}\Big[a & \left((t-x)_+^{h(v)-1/\alpha}
-(-x)_+^{h(v)-1/\alpha}\right) \nonumber \\
& +  b \left((t-x)_-^{h(v)-1/\alpha}
-(-x)_-^{h(v)-1/\alpha}\right) \Big]M(dx) \quad (t,v \in \bbbr).\label{rf}
\end{align}
Then $X(t,v)$ is well-defined since since for each $(t,v)$ 
the $\alpha$-th power of the integrand is 
Lebesgue integrable.  
For each fixed $v$ the process $X(\cdot,v)$ is just a 
linear stable fractional motion (\ref{lsfm}) so is $h(v)$-localisable,
with $X_{u}'(\cdot,v) = L_{\alpha,h(v)}(\cdot)$ for all
$u\in \bbbr$.  
Provided that 
$h(v)>1/\alpha$ it is in $C(\bbbr)$ and is strongly localisable.

\medskip
\begin{theo} (Linear multifractional stable motion)
Let  $U$ be  a closed interval with $u$ an interior point. Let
$0<\alpha<2$ and $h:U \to (0,1)$. Define
 $\{Y(t) : t\in U\}$  by $(\ref{lsmmdef})$ .

$(a)$ Assume that $h$ satisfies a $\eta$-H\"{o}lder condition at $u$
$$|h(v)-h(u)| \leq k|v-u|^{\eta}\quad (v \in U)$$ 
where 
$h(u) < \eta \leq 1$. Then 
$Y$ is $h(u)$-localisable at $u$ with local form 
$Y'_u=L_{\alpha,h(u)}$.

$(b)$ If $1<\alpha<2$ and $h$ is differentiable with 
$1/\alpha<h(u)<1$ and   
\begin{equation}
    |h'(v)-h'(v')|
    \leq k|v-v'|^{\eta}
    \quad (v,v' \in U)
    \label{hol6}
\end{equation}
where $1/\alpha<\eta \leq 1$, then $Y$ is strongly $h(u)$-localisable at $u$
with local form 
$Y'_u=L_{\alpha,h(u)}$.
\end{theo}

\noindent{\it Proof.}  
For brevity of exposition we give the proof in the case of
well-balanced linear multifractional stable motion, that is with 
$a=b=1$ in (\ref{lsmmdef}) and (\ref{rf}); the general
case is very similar. Thus we take
$$f(t,v,x) = |t-x|^{h(v)-1/\alpha}-|x|^{h(v)-1/\alpha}$$
in Theorem 7.3 (when $h(v)=1/\alpha$ such expressions are 
interpreted as $\1_{[0,t]}(x)$ where $\1_{[0,t]}$ is an indicator 
function). Then $X(t,v) = \int f(t,v,x) M(dx)$ and 
$Y(t) = \int f(t,t,x) M(dx)$.

(a)
By continuity, we may assume that $U$ is a sufficiently small interval 
to ensure that $h(v)< \eta$ for all $v \in U$. 
Fix $h_{-},h_{+}$ such that 
$0<  h_{-} < h(v) < h_{+}<1$ for all $v \in U$. 
Then for each $t,v,v',x \in U$ with $x \neq 0,x \neq t$, 
the mean value theorem gives
\begin{eqnarray}
|f(t,v,x)&-&f(t,u,x)| \nonumber\\
& =& \left| |t-x|^{h(\cdot)-1/\alpha}\log|t-x| 
-|x|^{h(\cdot)-1/\alpha}\log|x|\right||h(v)-h(u)|\nonumber\\
& \leq & \left| |t-x|^{h(\cdot)-1/\alpha}\log|t-x| 
-|x|^{h(\cdot)-1/\alpha}\log|x|\right|k|v-u|^{\eta},\label{est1}
\end{eqnarray}
where $h(\cdot) \equiv h(t,v,x) \in [h(v),h(u)]$.  But
$$k\left| |t-x|^{h(\cdot)-1/\alpha}\log|t-x| 
-|x|^{h(.)-1/\alpha}\log|x|\right| \leq k_{1}(t,x)$$
for all $t \in U, x\in \bbbr$, where
\begin{equation}
k_{1}(t,x) = \left\{
\begin{array}{ll}
    c_{1}\max\left\{1, |t-x|^{h_{-}-1/\alpha} 
    +|x|^{h_{-}-1/\alpha}\right\}
    & ( |x| \leq 1 + 2 \max_{t \in U}|t|) \\
    c_{2}|x|^{h_{+}-1/\alpha-1} 
     & ( |x| > 1 + 2 \max_{t \in U}|t|) 
\end{array}
\right.\label{ktx}
\end{equation}
for appropriately chosen constants $c_{1}$ and $c_{2}$. Then 
$\int k_{1}(t,x)^{\alpha}dx$ is finite and uniformly bounded for $t\in U$,
so as $X(\cdot,u)$ is $h(u)$-localisable at $u$, Theorem 7.3(a) gives
that $Y =\{ X(t,t): t\in U\}$ is $h(u)$-localisable at $u$ with
local form $Y'_u(\cdot)=X'_u(\cdot,u) = L_{\alpha,h(u)}(\cdot)$.

(b) We may assume that $U$ is small enough and $h_{-},h_{+}$ are 
chosen so that 
$0<  1/\alpha< h_{-} < h(v) < h_{+}<1$ for all $v \in U$.
A similar estimate to (\ref{est1}) on the derivatives gives 
\begin{eqnarray}
\left|f_{v}(t,v,x)
-f_{v}(t,v',x)\right|
& \leq & \left[|h'(v)||h'(v(\cdot))||v-v'| + |h(v)-h(v')|\right]k_{1}(t,x)\nonumber\\
& \leq & c_{1} k_{1}(t,x) |v-v'|^{\eta},\label{est2}
\end{eqnarray}
for $t,v,v' \in U, x\in \bbbr$, 
where $k_{1}(t,x)$ is as in (\ref{ktx}), so (\ref{hol5}) is satisfied.
Moreover,
\begin{eqnarray}
|f_{v}(t,v,x)
&-&f_{v}(t',v,x)| \nonumber\\
& =& |h'(v)|\big| |t-x|^{h(v)-1/\alpha}\log|t-x| 
-|t'-x|^{h(v)-1/\alpha}\log|t'-x|\big|\nonumber\\
& \leq & k_{2}(t,t',x),\label{est3}
\end{eqnarray}
for $t,t',v \in U, x\in \bbbr$, where
\begin{equation}
k_{2}(t,t',x) = \left\{
\begin{array}{ll}
    c_{2}|t-t'|^{h_{-}-1/\alpha}
    & ( |x-\textstyle{\frac{1}{2}}(t-t')| \leq |t-t'|) \\
    c_{3}|x-\textstyle{\frac{1}{2}}(t-t')|^{h_{+}-1/\alpha-1}|t-t'| 
     & ( |x-\textstyle{\frac{1}{2}}(t-t')| > |t-t'|) 
\end{array}
\right. \label{k2tx}
\end{equation}
for constants $c_{2},c_{3}$.
Then $\| k_{2}(t,t',\cdot)\|_{\alpha} \leq c_{4}|t-t'|^{h_{-}}$, so 
(\ref{hol7}) is satisfied taking $\eta =  h_{-}$.  Strong localisability
 follows from Theorem 7.3(c).
\Box
\medskip

To conclude this section we examine stationary moving average
processes.  These provide examples of localisable $\alpha$-stable 
processes of a rather
different nature being  stationary processes and not based on 
existing sssi processes.

\begin{prop}
Let $0<\alpha\leq 2$, let $g \in {\cal F}_{\alpha}$ and let $M$ be a
symmetric $\alpha$-stable
measure on $\bbbr$ with control measure ${\cal L}$.
Define the stationary process $Y$ by   
\begin{equation}
Y(t) = \int g(t-x)M(dx) \quad (t \in \bbbr).\label{ma}
\end{equation}
Suppose that there exist jointly measurable functions $h(t,.) \in 
{\cal F}_{\alpha}$ such that 
\begin{equation}
\lim_{r\to 0} \int\left|\frac{g(r(t+z))-g(rz)}{r^{\gamma}} -
h(t,z)\right|^{\alpha}dz=0 \label{malim}
\end{equation}
 for all $t\in \bbbr$, where $\gamma+(1/\alpha) >0$.  Then $Y$ is
$(\gamma+(1/\alpha))$-localisable at all $u\in \bbbr$ with local form
$Y_{u}= \{\int h(t,z)M(dz): t\in \bbbr\}$.  
\end{prop}

\noindent{\it Proof.} Using stationarity followed by a change of
variable $z=-x/r$ and the self-similarity of $M$,
\begin{eqnarray*}
Y(u+rt)-Y(u) &=& Y(rt)-Y(0) \\
&=& \int (g(rt-x)-g(-x))M(dx)\\
&=& r^{1/\alpha}\int (g(r(t+z))-g(rz))M(dz)
\end{eqnarray*}
where equality is in finite dimensional distributions.
Thus 
$$\frac{Y(u+rt)-Y(u)}{r^{\gamma+1/\alpha}} -\int h(t,z)M(dz)
= \int \left(\frac{g(r(t+z))-g(rz)}{r^{\gamma}} -
h(t,z)\right)M(dz).$$
By \cite[Proposition 3.5.1]{ST} and (\ref{malim}),
$r^{-\gamma-1/\alpha}(Y(u+rt)-Y(u))\to 
\int h(t,z)M(dz)$ in probability and thus in finite dimensional
distributions.  
\Box

\medskip

A particular instance of (\ref{ma}) is the reverse Ornstein-Uhlenbeck
process, see \cite[Section 3.6]{ST}. 

\begin{theo} (Reverse Ornstein-Uhlenbeck
process)
Let $\lambda>0$ and $0<\alpha \leq 2$ and let  $M$ be an $\alpha$-stable
measure on $\bbbr$ with control measure ${\cal L}$.
The stationary process defined by    
\begin{equation*}
Y(t) = \int_{t}^{\infty} \exp(-\lambda (x-t))M(dx) \quad (t \in \bbbr)
\end{equation*}
has a version in $D(\bbbr)$ that is $1/\alpha$-localisable
at all $u \in \bbbr$ with $Y_{u}' =L_{\alpha}$, where 
$L_{\alpha}$ is $\alpha$-stable L\'{e}vy motion.
\end{theo}

\noindent{\it Proof.} 
The process $Y$ is a stationary Markov process which has a version
in $D(\bbbr)$ see \cite[Remark 17.3]{Sa}.  It is a moving average process 
taking $g(x) = \exp(\lambda x) \1_{(-\infty,0]}(x)$ in (\ref{ma}).  It is
easily verified using the dominated convergence theorem 
that $g$ satisfies (\ref{malim}) with $\gamma = 0$ and
$h(t,z) = -\1_{[-t,0]}(z)$, so Proposition 7.5 gives the conclusion
with  $Y_{u}'(t) = -M ([-t,0]) = L_{\alpha}(t)$.  
\Box


\section{Sums over Poisson processes}
\setcounter{equation}{0}
\setcounter{theo}{0}
\medskip

In the next section we will set up `multistable processes', that is
$\alpha$-stable processes where $\alpha$ is allowed to vary with $t$.
For this it is convenient to express the random field $X(t,v)$ as a sum over 
a suitable Poisson point process.

In this section we bring together the basic properties of Poisson sums
that we need.  Let $(E,{\cal E},m)$ be a $\sigma$-finite measure 
space.  We work
throughout with a Poisson point process $\Pi$ on 
$E \times \bbbr$, with
 mean measure $m \times {\cal L}$ where ${\cal L}$ is 
Lebesgue measure.  Thus 
$\Pi$ is a random countable subset of $E \times \bbbr$ such that, writing
$N(A)$ for the
number of points in a measurable $A \subset E \times \bbbr$, the random
variable
$N(A)$ has a Poisson distribution of mean $(m \times {\cal L})(A)$
with
$N(A_{1}),\ldots,N(A_{n})$ independent
for disjoint $A_{1},\ldots,A_{n}\subset \bbbr^{2}$, see \cite{Ki}. 

We define a quasinorm on certain spaces of measurable functions on $E$. 
For $ 0<a \leq b <2$ let
$${\cal F}_{a,b}\equiv {\cal F}_{a,b}(E,{\cal E},m)
= \{f: f \mbox{ is $m$-measurable with } \|f\|_{a,b}<\infty\}$$
where
\begin{equation}
\|f\|_{a,b} = \left(\int_{E}|f(x)|^{a}m(dx)\right)^{1/a}
+  \left(\int_{E}|f(x)|^{b}m(dx)\right)^{1/b}.\label{norm}
\end{equation}
(Of course $\|\,\|_{a,b}$ is a norm if $1 \leq a \leq b$.) Note that
if $a \leq a'\leq b' \leq b$ then ${\cal F}_{a,b} \subset {\cal
F}_{a',b'}$ and $\|f\|_{a',b'} \leq c \|f\|_{a,b}$ where $c$ depends
on $a,a',b',b$. Moreover, 
${\cal F}_{a,a} ={\cal F}_{a} $.

The following estimate will be useful. Note that 
expressions such as (\ref{gest1}) have two parts since we need to control the growth 
of $g(x,y)$ at both small and large values of $y$.

\begin{lem}
Let  $g: E \times \bbbr \to \bbbr$ be ${\cal L}^{2}$-measurable and suppose
that
\begin{equation}
|g(x,y)| \leq h(x)\left(|y|^{-1/a} + |y|^{-1/b}\right) \label{gest1}
\end{equation}
where $h \in {\cal F}_{a,b}$ for some $0<a \leq b<2$.
Then there is a constant $c$ depending only on $a$ and $b$ such that
\begin{equation}
\int\int \sin^{2}({\textstyle\frac{1}{2}}\th g(x,y))m(dx)dy
\leq c\left(\th^{a}\int|h(x)|^{a}m(dx) +
\th^{b}\int|h(x)|^{b}m(dx)\right)\quad (\th \geq 0). \label{sinest}
\end{equation}    
\end{lem}    

\noindent {\it Proof.}
We have
\begin{eqnarray}
\int\int \sin^{2}({\textstyle\frac{1}{2}}\th g(x,y))m(dx)dy 
& \leq & \int\int \min
\left\{{\textstyle\frac{1}{4}}\th^{2}|g(x,y)|^{2}, 1\right\}m(dx) dy
    \nonumber\\
     & \leq & 
     c_{1}\int\int \min\left\{\th^{2}|h(x)|^{2}|y|^{-2/a},1\right\}m(dx) dy 
     \nonumber\\
     & & + c_{1}\int\int \min\left\{\th^{2}|h(x)|^{2}|y|^{-2/b},1\right\}m(dx)
     dy, \label{doub}  
\end{eqnarray}
 where $c_{1}$ is a constant, using (\ref{gest1}) and making a simple 
 estimate.  But
\begin{align*}
\int\int \min\big\{\th^{2}|h(x)|^{2} & |y|^{-2/a}, 1\big\} m(dx)dy \\
& \leq \int \left[\int_{|y| \leq |\th h(x)|^{a}} dy
      + \th^{2}|h(x)|^{2}\int_{|y| > |\th h(x)|^{a}}|y|^{-2/a} dy\right]
      m(dx)\\
& \leq  c_{2}\th^{a}\int|h(x)|^{a}m(dx)     
\end{align*}
where $c_{2}$ depends only on $a$, so along with a similar estimate with
$b$ replacing $a$,   (\ref{doub}) gives (\ref{sinest}).
\Box

\medskip

The next proposition gives criteria for the convergence of Poisson sums.
We write $(\X,\Y)$ for a random point of $E\times\bbbr$ of the Poisson
process $\Pi$.

\begin{prop}
Let  $g: E\times\bbbr \to \bbbr$ be $m \times {\cal L}$-measurable with 
\begin{equation}
|g(x,y)| \leq h(x)\left(|y|^{-1/a} + |y|^{-1/b}\right) \label{gest}
\end{equation}
where $h \in {\cal F}_{a,b}$.

$(a)$ If $0<a\leq b<1$ then the series
\begin{equation}
\s \equiv \sum_{(\X,\Y) \in \Pi} g(\X,\Y) \label{pois2}
\end{equation}
converges absolutely almost surely.

$(b)$ Suppose that  $0<a\leq b<2$ and that $g$ is
symmetric in the sense that
\begin{equation}
g(x,-y) = -g(x,y) \quad\quad (x,y) \in E\times\bbbr. \label{sym}
\end{equation}
Let $E_{n}$ be an increasing sequence of
$m$-measurable subsets of $E$ with $m(E_{n}) < \infty$ for all $n$
and $\cup_{n=1}^{\infty}E_{n}=E$ and write $R_{n}$ for the rectangle 
$\{(x,y)  : x \in E_{n},|y| \leq n\}\subset
E \times \bbbr$.  Then we may define
\begin{equation}
\s \equiv \sum_{(\X,\Y) \in \Pi} g(\X,\Y) 
= \lim_{n \to \infty} \sum_{(\X,\Y) \in \Pi \cap R_{n}} g(\X,\Y), \label{pois4}
\end{equation}
where the series converges almost surely.

$(c)$ Provided the symmetry condition $(\ref{sym})$ holds,
the characteristic function of $\s$, taking either definition $(\ref{pois2})$
or definition $(\ref{pois4})$, is
given by
\begin{equation}
\E(e^{i\th\s}) =
\exp\left(-2\int\int \sin^{2}({\textstyle\frac{1}{2}}\th g(x,y))m(dx)dy
\right) \quad (\th\in \bbbr).\label{cf}
\end{equation}
\end{prop}

\noindent {\it Proof.}
If $0<a \leq b <1$, (\ref{gest}) easily implies that 
$\int\min\{|g(x,y)|,1\}m(dx) dy <\infty$.  By Campbell's theorem
\cite[Section 3.2]{Ki} the random sum (\ref{pois2}) is absolutely
convergent almost surely with characteristic function
$$\E(e^{i\th\s}) =  
\exp\left(\int\int \left(e^{i\th g(x,y)}-1\right)m(dx)dy
\right) \quad (\th\in \bbbr).$$
If the symmetry condition (\ref{sym}) holds, this reduces to (\ref{cf}).

In case (b) 
write $\Sigma_{n}=\sum_{(\X,\Y) \in \Pi \cap R_{n}} g(\X,\Y) =
\sum_{(\X,\Y) \in \Pi} g(\X,\Y)\1_{R_{n}}(\X,\Y)$, where $\1_{R_{n}}$ is the
indicator function of  $R_{n}$.  Then by (\ref{gest})
$\int\min\{|g(x,y)\1_{R_{n}}(x,y)|,1\}m(dx) dy <\infty$, so
using Campbell's theorem just as before
\begin{eqnarray*}
\E(e^{i\th\Sigma_{n}})
& = &
\exp\left(-2\int\int \sin^{2}({\textstyle\frac{1}{2}}\th 
g(x,y))\1_{R_{n}}(x,y)m(dx)dy
\right)\\
& \to &
\exp\left(-2\int\int \sin^{2}({\textstyle\frac{1}{2}}\th g(x,y))m(dx)dy
\right),
\end{eqnarray*}
as $n \to \infty$ for all $\th$, by monotone convergence.
By (\ref{sinest}) there is a number $c_{1}>0$ such that 
$$1 \geq \exp\left(-2\int\int \sin^{2}({\textstyle\frac{1}{2}}\th g(x,y))m(dx)dy
\right) \geq  \exp(-2 c_{1}|\theta|^{a}) \geq 1 -
2c_{1}|\theta|^{a}$$
for $|\theta| \leq 1$, using that $1-e^{-x} \leq x$ if $x \geq 0$.
Thus $\lim_{n \to \infty} \E(e^{i\th\Sigma_{n}})$ exists for all $\th$
and is continuous at $\th=0$, so by L\'{e}vy's continuity theorem
\cite[Section 10.6]{Ei}, $\Sigma_{n}$ converges in distribution to a random
variable $\s$ with characteristic function (\ref{cf}).  

We may write
$$\lim_{n \to \infty} \sum_{(\X,\Y) \in \Pi \cap R_{n}} g(\X,\Y)
= \sum_{n=1}^{\infty}\,\sum_{(\X,\Y) \in \Pi \cap (R_{n}\setminus
R_{n-1})} g(\X,\Y)$$
(taking $R_{0}=\emptyset$), which is an infinite sum of independent random variables that converges in 
distribution, so by another theorem of L\'{e}vy \cite[Chapter 12]{Ei} it also
converges almost surely.
\Box

\medskip

\begin{prop}
Let $\s = \sum_{(\X,\Y) \in \Pi} g(\X,\Y)$ be as in $(\ref{pois2})$ or   
$(\ref{pois4})$ where $g(x,-y) = -g(x,y) $ and
\begin{equation}
|g(x,y)| \leq h(x)\left(|y|^{-1/a} + |y|^{-1/b}\right) \label{gest2}
\end{equation}
for some $h \in {\cal F}_{a,b}$.
  Then for $0<p< a$, 
\begin{equation}
\E|\s|^{p} \leq 
 c\|h\|_{a,b}^{p},\label{expzp}
\end{equation}
where $c$ depends only on $a,b$ and $p$.
\end{prop}

\noindent {\it Proof.}
A simple calculation using characteristic functions (see 
\cite[p.47]{Bil}) gives
\begin{eqnarray*}
    \P\{|\s| \geq \lambda\}& \leq &
    \frac{\lambda}{2}\int_{-2/\lambda}^{2/\lambda}\big(1-\E(\exp(i\theta
    \s))\big)d\theta 
    \\
    &= & 
    \frac{\lambda}{2}\int_{-2/\lambda}^{2/\lambda}
    \bigg(1-\exp\Big(-2\int\int \sin^{2}({\textstyle\frac{1}{2}}\th
    g(x,y))m(dx)dy\Big)\bigg) d\theta\\
     &\leq & 
    \frac{\lambda}{2}\int_{-2/\lambda}^{2/\lambda}
    \bigg(1-\exp\Big(-2c\big(\th^{a}\int|h(x)|^{a}m(dx) +
\th^{b}\int|h(x)|^{b}m(dx)\big)\Big)\bigg)d\theta\\
     &\leq & 
    c\lambda\int_{-2/\lambda}^{2/\lambda}
    \left(\th^{a}\int|h(x)|^{a}m(dx) +
\th^{b}\int|h(x)|^{b}m(dx)\right)d\theta \\
 &\leq & 
    c_{1} \lambda^{-a}\int|h(x)|^{a}m(dx) + c_{1}
    \lambda^{-b}\int|h(x)|^{b}m(dx)\\
 &\equiv &  \lambda^{-a}h_{a} + \lambda^{-b}h_{b}, 
\end{eqnarray*}
say, where $c_{1}$ depends only on $a$ and $b$, using (\ref{cf}) and (\ref{sinest}).
Then
\begin{eqnarray*}
\E|\s|^{p} &=& p\int_{0}^{\infty}\lambda^{p-1}\P(|\s| \geq \lambda) 
 d\lambda \\
 &\leq &  p\int_{0}^{\infty}  \lambda^{p-1}\min\{1, \lambda^{-a}h_{a}\}d\lambda 
 +  p\int_{0}^{\infty}  \lambda^{p-1}\min\{1, \lambda^{-b}h_{b}\}d\lambda\\
 &\leq &  p\int_{0}^{h_{a}^{1/a}}\lambda^{p-1}d\lambda 
    + p h_{a}\int_{h_{a}^{1/a}}^{\infty}\lambda^{p-a-1}d\lambda 
 + p\int_{0}^{h_{b}^{1/b}}\lambda^{p-1}d\lambda 
    +  ph_{b}\int_{h_{b}^{1/b}}^{\infty}\lambda^{p-b-1}d\lambda\\
    & \leq &  c_{2} (h_{a}^{p/a} + h_{b}^{p/b})\\
& \leq &  c\|h\|_{a,b}^{p}
\end{eqnarray*}    
where $c_{2},c$ depend on $a,b$ and $p$.
\Box

\medskip

We will sometimes need the following variant of Proposition 8.3.

\begin{cor}
Let $0<p <a <a_{1}<b_{1} < b <2$ and let $f_{1},f_{2} \in {\cal F}_{a,b}$. 
Let 
$$\s \equiv \sum_{(\X,\Y) \in \Pi} \left(
f_{1}(\X) \Y^{<-1/\alpha_{1}>} - f_{2}(\X) \Y^{<-1/\alpha_{2}>}\right),$$
where $a_{1}\leq \alpha_{1},\alpha_{2} \leq b_{1}$. Then 
$$\E|\s|^{p} \leq c\|f_{1}-f_{2}\|_{a,b}^{p}
+ c \|f_{2}\|_{a,b}^{p}|\alpha_{1}-\alpha_{2}|^{p}$$
where $c$ depends only on $a,a_{1},b,b_{1}$ and $p$.
\end{cor}

\noindent {\it Proof.}
Since
\begin{eqnarray*}
\s  & = &  \sum_{(\X,\Y) \in \Pi}(f_{1}(\X)  - f_{2}(\X)) \Y^{<-1/\alpha_{1}>} 
+f_{2}(\X) \sum_{(\X,\Y) \in \Pi} \left(\Y^{<-1/\alpha_{1}>} -  \Y^{<-1/\alpha_{2}>}\right) \\
     & = & \sum_{(\X,\Y) \in \Pi}(f_{1}(\X)  - f_{2}(\X)) \Y^{<-1/\alpha_{1}>} 
+f_{2}(\X) \sum_{(\X,\Y) \in \Pi}
(\alpha_{1}-\alpha_{2})\Y^{<-1/\alpha>}\alpha^{-2}\log |\Y|
\end{eqnarray*}
where $\alpha \in [\alpha_{1},\alpha_{2}]$, using the mean value
theorem,
the corollary follows from Proposition 8.3.
\Box
\medskip

Note that the introduction of $a_{1}$ and $b_{1}$ in Corollary 8.4 is 
necessitated by the `log' term to ensure uniformity of the constant $c$.

\section{Multistable processes}
\setcounter{equation}{0}
\setcounter{theo}{0}
\medskip

We now show how our approach may be used to construct multistable processes, 
that is
processes where the local stability index varies.
The development of this section
mirrors that of Section 7, but depends heavily on
the properties of Poisson sums derived in 
Section 8.
We seek an analogue of 
Theorem 7.3 but with the local form $Y_{u}'$ an $\alpha(u)$-stable process
with $\alpha(u)$ depending on $u$.  

We first define a random field
analogous to (\ref{field1}), but where the stable random measure 
$M$ is not allied to a
particular value of $\alpha$.  Whilst it would be possible to set up
a random measure that resembles an $\alpha(u)$-stable measure close to
$u$, this would be technically 
quite complicated.  We therefore favour an
alternative approach, using a representation by sums 
over Poisson processes.   In particular this permits $X(\cdot,v)$ to be specified 
using the same underlying Poisson process for different $v$.

As before $(E,{\cal E},m)$ is a 
$\sigma$-finite measure space, and 
$\Pi$ is a Poisson process on $E \times \bbbr$ with mean measure
 $m \times {\cal L}$.  
In the case of constant $\alpha$, with $M$  a symmetric $\alpha$-stable
random measure on $E$ with control measure $m$ and skewness $0$,
the stochastic integral 
(\ref{alint}) 
may be expressed as a Poisson process sum
\begin{equation}
I(f) = \int f(x) M(dx)\, = \, c(\alpha)
\sum_{(\X,\Y) \in \Pi} f(\X) \Y^{<-1/\alpha>}\quad (0<\alpha<2),\label{ppal}
\end{equation}
with the sum taken in the sense of (\ref{pois2}) or
(\ref{pois4}), and with 
\begin{equation}
    c(\alpha) =
\left(2\alpha^{-1}\Gamma(1-\alpha)
\cos(\textstyle{\frac{1}{2}}\pi
\alpha)\right)^{-1/\alpha},\label{calpha}
\end{equation}
    see
\cite[Section 3.12]{ST}. (As
before $a^{<b>} = \mbox{sign}(a)|a|^{b}$ and ${\cal L}$ is Lebesgue measure.)

Particularly relevant in  (9.1) is that the stability index $\alpha$
occurs only as an exponent of $\Y$, since the underlying Poisson process does
not depend on $\alpha$, so by varying this exponent we can vary the
stability index.
Thus the random field
\begin{equation}
X(t,v)= \sum_{(\X,\Y) \in \Pi} f(t,v,\X) \Y^{<-1/\alpha(v)>}\label{msfielda}
\end{equation}
gives rise to  a {\it multistable} process with varying $\alpha$,
of the form
\begin{equation}
Y(t) \equiv X(t,t)= \sum_{(\X,\Y) \in \Pi} f(t,t,\X)
\Y^{<-1/\alpha(t)>}.\label{msproc}
\end{equation}

 We first consider continuity and boundedness of the processes.

\begin{prop}
Let $U$ be a closed interval.
Let $X$ be the random field defined by
\begin{equation}
X(t,v)= \sum_{(\X,\Y) \in \Pi} f(t,v,\X) \Y^{<-1/\alpha(v)>}
\quad (t,v \in U)\label{msfield1}
\end{equation}
where $f(t,v,\cdot)\in {\cal F}_{a,b}$ are jointly measurable and
$\alpha: U \to (a,b)$ is continuous.

$(a)$  Suppose $0< a<\alpha(v)<b<1$ for $v \in U$.  If
\begin{equation}
\sup_{t,v \in U}|f(t,v,x)|\leq k(x), \label{supfin2}
\end{equation}
where $k \in {\cal F}_{a,b}$,
 then  $\{X(t,v): t,v \in U\}$ has a bounded version.
 
 If in addition $\{f(t,v,x) : x \in E\}$ is an equiuniformly
 continuous family
 for $t,v \in U$, then $X$ has a continuous version.
 
$(b)$ Suppose that $1< a<\alpha(v)<b<2$ for $v \in U$ 
and  $1/a<\eta\leq 1 $.  Suppose that
\begin{equation}
    |\alpha (v)-\alpha(v')| \leq k_{1} |v-v'|^{\eta}
    \quad (v,v' \in U),
    \label{holal}
\end{equation}
that
\begin{equation}
    \sup_{t,v \in U}\|f(t,v,\cdot)\|_{a,b}  < \infty,
    \label{bound}
\end{equation}
and 
\begin{equation}
    \|f(t,v,\cdot)-f(t',v',\cdot)\|_{a,b} 
    \leq k_{2}\left(|v-v'|^{\eta}+|t-t'|^{\eta}\right)
    \quad (t,t',v,v' \in U).
    \label{cty1s}
\end{equation}
Then $Y =\{ X(t,t): t\in U\}$ has a continuous version
satisfying an a.s. $\beta$-H\"{o}lder condition for all 
$ 0 <\beta <(\eta a-1)/a$.
\end{prop}

 \noindent{\it Proof.} 
 (a) From (\ref{msfield1}), for all $t,v \in U$,
\begin{eqnarray*}
    |X(t,v)| &\leq& \sum_{(\X,\Y) \in \Pi} |f(t,v,\X)|
    |\Y|^{-1/\alpha(v)}\\
&\leq & \sum_{(\X,\Y) \in \Pi} \sup_{t,v \in U}
|f(t,v,\X)|( |\Y|^{-1/a}+ |\Y|^{-1/b})\equiv Z
\end{eqnarray*}
where $Z$ is an a.s. finite random variable, by Proposition 8.2(a).
Thus $\{X(t,v): t,v \in U\}$ is a.s. bounded. 

Assuming also the equicontinuity condition, given $\epsilon>0$ we
may choose $r\geq 1$ such that $\|k(x)\1_{\{|x|> r\}}(x)\|_{a,b}<
\epsilon$, where $\1$ is the indicator function.   
By equiuniform continuity we may find $\delta>0$ such that for
all $x \in \bbbr$ and $|(t,v)-(t',v')|<\delta$ we have $|f(t,v,x)-f(t',v',x)|<
r^{-1/a}\epsilon$, and $|\alpha(v) - \alpha(v') |<\epsilon$.
Then if $|(t,v)-(t',v')|<\delta$,  making several estimates in the
obvious way,
\begin{eqnarray*}
    |X(t,v)&-&X(t',v')|
    \quad\leq \quad\sum_{(\X,\Y) \in \Pi} 
    |f(t,v,\X)\Y^{<-1/\alpha(v)>}-f(t',v',\X)\Y^{<-1/\alpha(v')>}| \\
   &\leq &\sum_{|\X| \leq r} 
    |f(t,v,\X)-f(t',v',\X)||\Y|^{-1/\alpha(v)}
    +2\sum_{|\X| > r} 
    \sup_{t,v \in U}|f(t,v,\X)||\Y|^{-1/\alpha(v)}\\
    & & + \sum_{(\X,\Y) \in \Pi} 
    |f(t',v',\X)||\Y^{<-1/\alpha(v)>}-\Y^{<-1/\alpha(v')>}|\\
 &\leq &\sum_{(\X,\Y) \in \Pi} 
    r^{-1/a}\epsilon \1_{\{|x|\leq r\}}(\X)|\Y|^{-1/\alpha(v)}
    +2\sum_{(\X,\Y) \in \Pi}k(\X)\1_{\{|x|> r\}}(\X) |\Y|^{-1/\alpha(v)}\\
    & & + \sum_{(\X,\Y) \in \Pi} 
    |f(t',v',\X)|\frac{1}{\alpha^{2}}|\Y|^{-1/\alpha}|\log|\Y||\,
    |\alpha(v)-\alpha(v')|\\    
  &\leq &\,\Big(\sum_{(\X,\Y) \in \Pi} 
    r^{-1/a}\epsilon \1_{\{|x|\leq r\}}(\X)
    +2\sum_{(\X,\Y) \in \Pi}k(\X)\1_{\{|x|> r\}}(\X) \\
    & & + \sum_{(\X,\Y) \in \Pi} 
    c_{1}|k(\X)|\frac{1}{a^{2}}
    \epsilon\Big)(|\Y|^{-1/a}+|\Y|^{-1/b}) \equiv Z_{\epsilon}
\end{eqnarray*}
where  
$Z_{\epsilon}$ is a random variable, and we have used the mean
value theorem in the third term of the sum with $\alpha \in
[\alpha(v),\alpha(v')]$.  Fix $0<p<\alpha$.  
By (\ref{expzp}) there is a constant 
$c$ independent of $\epsilon$ such that
$$\E|Z_{\epsilon}|^{p} \leq c \epsilon^{p}.$$
The proof is completed just as in the proof of Proposition 7.1(a).
 
 (b) We estimate
\begin{equation}
X(t,v)-X(t',v') = (X(t,v)-X(t,v')) + (X(t,v')-X(t',v'))
\quad (t,t',v,v' \in U) \label{diffxx}
\end{equation}  
by considering its two parts in turn. Firstly
$$ X(t,v)-X(t,v')
 =   \sum_{(\X,\Y)\in \Pi} 
\left(f(t,v,\X) \Y^{<-1/\alpha(v)>}- f(t,v',\X) \Y^{<-1/\alpha(v')>}
\right)$$
Thus  Corollary 8.4 gives, for $0<p<a$,
\begin{eqnarray}
   \E |X(t,v)-X(t,v')|^{p}  
     &\leq& c_{1}  \|f(t,v,\cdot)-f(t,v',\cdot)\|_{a,b}^{p}
     + c_{1}  \|f(t,v',\cdot)\|_{a,b}^{p} 
     |\alpha(v)-\alpha(v')|^{ p} \nonumber\\
     &\leq&  c_{2} |v-v'|^{\eta p}
         \label{1stctys}
\end{eqnarray}
by (\ref{cty1s}), (\ref{bound}) and (\ref{holal}).

For the second term of (\ref{diffxx})
$$ X(t,v)-X(t',v)
=    \sum_{(\X,\Y)\in \Pi} 
\left(f(t,v,\X) - f(t',v,\X) \right)\Y^{<-1/\alpha(v)>}.
$$
Then 
$$|(f(t,v,x) - f(t',v,x) )y^{<-1/\alpha(v)>} |
\leq |f(t,v,x) - f(t',v,x) |\left(|y|^{-1/a}+|y|^{-1/b}\right)$$
so, for $0<p<a$,  Proposition 8.3 and (\ref{cty1s}) give
\begin{eqnarray*}
   \E |X(t,v)- X(t',v)|^{p}
   & \leq & c_{3}\|f(t,v,\cdot) - f(t',v,\cdot) \|_{a,b}^{p} \\
& \leq & c_{4}|t-t'|^{\eta p}.
\end{eqnarray*}

Combining with (\ref{1stctys}) we estimate (\ref{diffxx}) 
to get, for $t,t',v,v' \in U$,  
\begin{equation*}
   \E |X(t,v)- X(t',v')|^{p}
   \leq c_{5} (|v-v'|^{\eta p}+|t-t'|^{\eta p}). 
\end{equation*} 
Specialising, 
\begin{equation*}
   \E |Y(t)- Y(t')|^{p}=\E |X(t,t)- X(t',t')|^{p}  
   \leq 2c_{5}|t-t'|^{\eta p} 
\end{equation*} 
for $t,t' \in U$.  

Since $\eta > 1/a$ we may choose $0<p<a$ 
such that $\eta p >1$.   Kolmogorov's Theorem 5.2 gives that 
$\{Y(t): t\in U\}$ has a continuous version that is a.s.
$\beta$-H\"{o}lder for all 
$ 0 <\beta <(\eta p-1)/p$ for all $p<a$. 
 \Box
\medskip
 
We come to the main result on the localisability of processes with
varying stability index.

\begin{theo}
Let $U$ be a closed interval with $u$ an interior point and let $0<a < b <2$.  
Let $X$ be the random field defined by
\begin{equation}
X(t,v)= \sum_{(\X,\Y) \in \Pi} f(t,v,\X) \Y^{<-1/\alpha(v)>}
\quad (t,v \in U)\label{msfield}
\end{equation}
where $f(t,v,\cdot)\in {\cal F}_{a,b}$ are jointly measurable and
$\alpha: U \to (a,b)$.
   
$(a)$ Suppose  $X(\cdot,u)$ 
is $h$-localisable at $u$ for  $h>0$.
Suppose that
$\sup_{t\in U}
\|f(t,u,\cdot)\|_{a,b}<\infty$, and that for some $\eta >h$ 
\begin{equation}
    |\alpha(v) - \alpha(u) | \leq k_{1} |v-u|^{\eta}
    \quad (v \in U), \label{mshol2}
    \end{equation}
and
\begin{equation}
    \|f(t,v,\cdot) - f(t,u,\cdot)\|_{a,b} \leq k_{2} |v-u|^{\eta}
    \quad (t,v \in U). \label{hol2ms}
    \end{equation}
Then $Y =\{ X(t,t): t\in U\}$ is $h$-localisable at $u$ with
local form $Y'_u(\cdot)=X'_u(\cdot,u).$

$(b)$ Suppose that  $0<\alpha(u)<1$ and that  $X(\cdot,u)$ is 
strongly $h$-localisable in $C(\bbbr)$ (resp. $D(\bbbr)$)
at $u$.  Suppose that for some $\eta >h$ 
\begin{equation}
    |\alpha(v) - \alpha(u) | \leq k_{1} |v-u|^{\eta}
    \quad (v \in U), \label{mshol2b}
    \end{equation}
 and 
\begin{equation}
    |f(t,u,x)| \leq k_{2}(x)
    \quad (t \in U, x \in E),\label{hol4a}
\end{equation} 
and 
    \begin{equation}
    |f(t,v,x) - f(t,u,x)| \leq k_{3}(x) |v-u|^{\eta}
    \quad (t,v \in U, x \in E),\label{hol4ms}
\end{equation}
 where $k_{2}(\cdot),k_{3}(\cdot)\in{\cal F}_{a,b}$.
 If $Y =\{ X(t,t): t\in U\}$ has a version in $C(U)$ (resp. $D(U)$)
 then
 $Y$ is strongly $h$-localisable 
 in $C(\bbbr)$ (resp. $D(\bbbr)$)  at $u$ with
$Y'_u(\cdot)=X'_u(\cdot,u).$

$(c)$ Suppose that  $1<\alpha(u)<2$ and  that $X(\cdot,u)$ 
is strongly $h$-localisable
 in $C(\bbbr)$ (resp. $D(\bbbr)$)  at $u$. Let $\eta$ satisfy
 $1/\alpha(u) < \eta \leq 1$.  Suppose 
that $\alpha$ is continuously differentiable on $U$ with
\begin{equation}
  |\alpha'(v)-\alpha'(v')|\leq k_{1}|v-v'|^{\eta}
    \quad (v,v' \in U). \label{alhol}
\end{equation}
Suppose that 
the partial derivatives $f_{v}(t,v,\cdot) \in {\cal F}_{a,b}$ 
for all $t,v \in U$, and the following estimates hold: 
\begin{equation}
   \sup_{t \in U}\|f(t,u,\cdot)\|_{a,b}< \infty, \quad
   \sup_{t,v \in U}\|f_{v}(t,v,\cdot)\|_{a,b} < \infty,
    \label{holbms}
\end{equation}
\begin{equation}
   \| f(t,v,\cdot)- f(t',v,\cdot)\|_{a,b}
    \leq k_{2}|t-t'|^{\eta}
    \quad (t,t',v \in U),
    \label{new}
\end{equation}
\begin{equation}
   \left|f_{v}(t,v,x)
    - f_{v}(t,v',x)\right| 
    \leq k_{3}(t,x) |v-v'|^{\eta}
    \quad (t,v,v' \in U, x \in E),
    \label{holcms}
\end{equation}
and
\begin{equation}
    \left|f_{v}(t,v,x)
    - f_{v}(t',v,x)\right| 
    \leq k_{4}(t,t',x)
    \quad (t,t',v \in U, x \in E),
    \label{holdms}
\end{equation}where 
$\sup_{t \in U}\| k_{3}(t,\cdot)\|_{a,b}< \infty$, and
$\| k_{4}(t,t',\cdot)\|_{a,b}\leq k|t-t'|^{\eta}$ 
for all $t,t' \in U$.
Then $Y =\{ X(t,t): t\in U\}$ is strongly $h$-localisable 
in $C(\bbbr)$ at $u$ with
$Y'_u(\cdot)=X'_u(\cdot,u).$
\end{theo}

\medskip
\noindent{\it Proof.}
(a) We have
\begin{equation}
    X(t,v)-X(t,u) = \sum_{(\X,\Y) \in \Pi}\left(f(t,v,\X)\Y^{<-1/\alpha(v)>}
    - f(t,u,\X) \Y^{<-1/\alpha(u)>}\right).
 \label{msdiff}
\end{equation}
With $0<p<a$,  Corollary 8.4 gives that there are constants
$c_{1},c_{2}$  such that
\begin{eqnarray*}
\E |X(t,v)-X(t,u)|^{p} &\leq & 
c_{1} \| f(t,v,\cdot)-f(t,u,\cdot)\|_{a,b}^{p} 
+ c_{1} |\alpha(v)-\alpha(u)|^{p}\| f(t,u,\cdot)\|_{a,b}^{p} \\
&\leq & c_{2}|v-u|^{\eta p}, 
\end{eqnarray*}
for all $t,v \in U$, by (\ref{mshol2}) and (\ref{hol2ms}).
Part (a) now follows from
Theorem 4.1. 

(b) We may assume that $a< \alpha(v) < b <1$ for $v \in U$,
if necessary using the continuity of $\alpha$ to replace $U$ by a 
subinterval to decrease the value of $b$.
Splitting and estimating (\ref{msdiff}), using the mean value theorem
as in the proof of Corollary 8.4, we
get, with  $\alpha'\in [\alpha(v),\alpha(u)]$ (where $\alpha'$ depends on
$v$),
\begin{eqnarray*}
    |X(t,v)-X(t,u)|&\leq &\sum_{(\X,\Y) \in \Pi}|f(t,v,\X) - f(t,u,\X)|
    |\Y|^{-1/\alpha(v)} \\
&&+\sum_{(\X,\Y) \in \Pi}|f(t,u,\X)| 
|\alpha(v)-\alpha(u)||\Y|^{-1/\alpha'}\alpha'^{-2}|\log
|\Y|\,| \\
&\leq &|v-u|^{\eta}\sum_{(\X,\Y) \in
\Pi}|k_{3}(\X)|(|\Y|^{-1/a}+ |\Y|^{-1/b})\\ 
&&+c_{1}|v-u|^{\eta}\sum_{(\X,\Y) \in \Pi}|k_{2}(\X)| 
(|\Y|^{-1/a}+ |\Y|^{-1/b}) 
\end{eqnarray*}
for all $t,v \in U$, using (\ref{mshol2b})-(\ref{hol4ms}). By Proposition 8.2(a) 
$|k_{2}(\X)|(|\Y|^{-1/a}+ |\Y|^{-1/b})$ and $|k_{3}(\X)|(|\Y|^{-1/a}+
|\Y|^{-1/b})$ are a.s. finite random variables, so
(\ref{conda}) holds and Theorem 5.1 implies that
 $Y$ is strongly localisable at $u$.
 
(c) The conditions of Proposition 9.1(b) are easily checked, so $Y$
has a continuous version.  
We may assume that $1< a< \alpha(v) < b$ for $v \in U$ and that $1/a
<\eta$, using continuity of $\alpha$ 
to replace $U$ by a 
subinterval and to increase the value of $a$ if necesssary.

Define
$$Z(t,v) =\frac{X(t,v)-X(t,u)}{v-u} \quad (t,v \in U, v\neq u);$$
again we use Kolmogorov's criterion to
show that $\{Z(v,v): v \in U\}$ is almost surely bounded to get
(\ref{conda}).  
We write
\begin{equation}
Z(t,v) = Z_{1}(t,v) + Z_{2}(t,v)
\quad (t,v \in U) \label{zsum}
\end{equation} 
where
\begin{equation}
Z_{1}(t,v)   =\sum_{(\X,\Y)\in \Pi}g(t,v,\X)\Y^{<-1/\alpha(v)>}\,\,
\mbox{ with }\,\, g(t,v,x) = \frac{f(t,v,x)-f(t,u,x)}{v-u} \label{Z1}
\end{equation}
and
\begin{equation}
Z_{2}(t,v)  =\sum_{(\X,\Y)\in \Pi}f(t,u,\X)
\frac{\Y^{<-1/\alpha(v)>}-\Y^{<-1/\alpha(u)>}}{v-u}.
\label{Z2}
\end{equation}
For $p<a$ we estimate $\E|Z(t,v)-Z(t',v')|^{p}$ by breaking it into four parts.

(i)  Applying Lemma 7.2 to (\ref{holcms}) gives
$|g(t,v,x)-g(t,v',x)| \leq 2^{\eta}k_{3}(t,x)|v-v'|^{\eta}$. Thus
Corollary 8.4 on  (\ref{Z1}) and then (\ref{holbms}) 
with the mean value theorem gives
\begin{eqnarray}
\E|Z_{1}(t,v)-Z_{1}(t,v')|^{p}
&\leq& c_{1} \|g(t,v,\cdot)-g(t,v',\cdot)\|_{a,b}^{p} 
+ c_{1} \|g(t,v',\cdot)\|_{a,b}^{p} |\alpha(v)- \alpha(v')|^{p}
\nonumber\\
&\leq& c_{2} \|k_{3}(t,\cdot)\|_{a,b}^{p}|v-v'|^{\eta p} 
+ c_{2} \sup_{t,v \in U}\|f_{v}(t,v,\cdot)\|_{a,b}^{p} |v-v'|^{\eta p}
\nonumber\\
&\leq& c_{3}|v-v'|^{\eta p}.
\label{E1}
\end{eqnarray}

(ii) Using the mean value theorem and (\ref{holdms})
\begin{eqnarray*}
|g(t,v,x)-g(t',v,x)| &=& \frac{1}{|v-u|}
\left|(f(t,v,x)- f(t',v,x))- (f(t,u,x)-f(t',u,x))\right|\\
&=&\left|f_{v}(t,v_{1},x)- f_{v}(t',v_{1},x)\right|\\
&\leq & k_{4}(t,t',x)
\end{eqnarray*}
where $v_{1}\in (v,u)$ depends on $t,t',v$ and $x$.   Proposition 8.3 with (\ref{Z1}) now gives
\begin{equation}
\E|Z_{1}(t,v)-Z_{1}(t',v)|^{p}
\leq c_{4} \|k_{4}(t,t',\cdot)\|_{a,b}^{p} \leq c_{5}|t-t'|^{\eta p}. \label{E2}
\end{equation}

(iii)  Turning to (\ref{Z2}), a simple estimate using (\ref{alhol})
gives that
$$\left|\left[\frac{d}{dv}y^{<-1/\alpha(v)>}\right]_{v'}^{v}\right|
\leq c_{6}|v-v'|^{\eta}\left(|y|^{-1/a}+|y|^{-1/b}\right).$$
By Lemma 7.2
$$\left|\frac{y^{<-1/\alpha(v)>}-y^{<-1/\alpha(u)>}}{v-u}
-\frac{y^{<-1/\alpha(v')>}-y^{<-1/\alpha(u)>}}{v'-u}
\right| \leq 2^{\eta}c_{6}|v-v'|^{\eta}\left(|y|^{-1/a}+|y|^{-1/b}\right).$$
Thus Proposition 8.3 applied to (\ref{Z2}) gives
\begin{equation}
\E|Z_{2}(t,v)-Z_{2}(t,v')|^{p}
\leq  c_{7} \|f(t,u,\cdot)\|_{a,b}^{p}|v-v'|^{\eta p}
\leq c_{8}|v-v'|^{\eta p}.
\label{E3}
\end{equation}

(iv)  Finally, a further mean value estimate gives
\begin{align*}
\Big|(f(t,v,x)-
f(t',v,x))&\frac{y^{<-1/\alpha(v)>}-y^{<-1/\alpha(u)>}}{v-u}\Big|\\ 
&\leq \,\, c_{9}|(f(t,v,x)-f(t',v,x))|\left(|y|^{1/a}+|y|^{1/b}\right).
\end{align*}
By Proposition 8.3 and (\ref{new})
\begin{equation}
\E|Z_{2}(t,v)-Z_{2}(t',v)|^{p}
\leq c_{10} \|f(t,v,\cdot)-f(t',v,\cdot)\|_{a,b}^{p} \leq c_{11}|t-t'|^{\eta p}. \label{E4}
\end{equation}

Taking  (\ref{zsum}) with (\ref{E1}), (\ref{E2}), 
(\ref{E3}) and (\ref{E4}), we
conclude that for some $c_{12}$ independent of $t,t',v,v' \in U$,  
\begin{equation}
   \E |Z(t,v)- Z(t',v')|^{p}
   \leq c_{12} (|v-v'|^{\eta p}+|t-t'|^{\eta p}) \label{bothms}  
\end{equation} 
if $0<p<a$. Specialising, 
\begin{equation}
   \E |Z(v,v)- Z(v',v')|^{p}  
   \leq 2c_{12}|v-v'|^{\eta p} \label{both1ms}  
\end{equation} 
for $v,v' \in U$.  

Since $\eta > 1/a$ we may choose $0<p<a$ 
such that $\eta p >1$.  Applying  Kolomogorov's Theorem 5.2 to
$\{Z(v,v) :v \in U\}$ we conclude that 
$ {\displaystyle Z(v,v) =\frac{X(v,v) - X(v,u)}{v-u}}$ has a version that is 
a.s. bounded.  Thus (\ref{conda}) holds and strong 
localisability follows from Theorem 5.1.
 \Box
\medskip

We now show how Theorem 9.2 may be used to construct some specific
multistable processes.  In these examples we take $(E,{\cal E},m)$ 
to be Lebesgue measure on $\bbbr$ so that
$\Pi$ is the Poisson process on $\bbbr^{2}$ with mean measure
 ${\cal L}^{2}$.  
We first construct a multistable analogue of the
linear multifractional motion of Theorem 7.4.

\begin{theo} {\it (Linear multistable multifractional motion)}.
Let  $a : \bbbr \to \bbbr^{+}$
be $C^{1}$ and $\alpha:\bbbr \to (0,2)$ and
$h:\bbbr \to (0,1)$ be $C^{2}$. Define
\begin{equation} 
Y(t)  = a(t)c(\alpha(t))\sum_{(\X,\Y)\in\Pi}\Y^{<-1/\alpha(t)>}
\left(|t-\X|^{h(t)-1/\alpha(t)}-|\X|^{h(t)-1/\alpha(t)}\right) \quad
(t\in \bbbr).
\label{lmm}
\end{equation}

$(a)$ The process  $Y$ is $h(u)$-localisable
at all $u \in \bbbr$, with $Y_{u}' = a(u)L_{\alpha(u),h(u)}$,
where 
$L_{\alpha,h}$ is linear stable motion.

$(b)$ If $u$ is such that $h(u)> 1/\alpha(u)$ then $Y$ 
is strongly $h(u)$-localisable
in $C(\bbbr)$ at $u$, with $Y_{u}' = a(u)L_{\alpha(u),h(u)}$.
\end{theo}

\noindent{\it Proof.}  By the amplitude result, Proposition 3.2, the term  
$a(t)c(\alpha(t))$ in (\ref{lmm}) does not affect localisability, so 
it is enough to prove the result with $a(t)c(\alpha(t))=1$.
Define a random field by
\begin{eqnarray*}
X(t,v) &=&\sum_{(\X,\Y)\in\Pi}\Y^{<-1/\alpha(v)>}
\left(|t-\X|^{h(v)-1/\alpha(v)}-|\X|^{h(v)-1/\alpha(v)}\right)  \quad 
(t,v \in \bbbr)\\
&=& \sum_{(\X,\Y)\in\Pi} f(t,v,\X) \Y^{<-1/\alpha(v)>},
\end{eqnarray*}
where 
$$f(t,v,x) = 
\left(|t-x|^{h(v)-1/\alpha(v)}-|x|^{h(v)-1/\alpha(v)}\right).$$
Then 
\begin{align*}
    f_{v} & (t,v,x) \\
& =\left(|t-x|^{h(v)-1/\alpha(v)}\log
|t-x|-|x|^{h(v)-1/\alpha(v)}\log
|x|\right)\left(h'(v)+\alpha'(v)/\alpha(v)^{2}\right).
\end{align*}
Given $u \in \bbbr$ we may use continuity of $h$ and $\alpha$
to choose $U$ to be a small enough 
closed interval with $u$ an
interior point, and numbers $a,b,h_{-},h_{+}$, such that
$0<a <\alpha(v) < b<2 $ and  $0<h_{-} < h(v) < h_{+} <1$
for all $v \in U$, and such that
$\frac{1}{a}-\frac{1}{b}< h_{-}<h_{+}<1-(\frac{1}{a}-\frac{1}{b})$.
A similar argument to that of 
Theorem 7.4 gives that 
\begin{equation}
|f(t,v,x)|,
|f_{v}(t,v,x)| 
\leq  k_{1}(t,x)\quad (t,v\in U, x \in \bbbr)
\end{equation}
and
\begin{equation}
|f(t,v,x)-f(t,v',x)|,
|f_{v}(t,v,x)
-f_{v}(t,v',x)| \\
\leq k_{1}(t,x) |v-v'|\, (t,v,v'\in U, x \in \bbbr) \label{fdiff}
\end{equation}
 where
\begin{equation}
k_{1}(t,x) = \left\{
\begin{array}{ll}
    c_{1}\max\{1, |t-x|^{h_{-}-1/a} +|x|^{h_{-}-1/a}\}
    & ( |x| \leq 1 + 2 \max_{t \in U}|t|) \\
    c_{2}|x|^{h_{+}-1/b-1} 
     & ( |x| > 1 + 2 \max_{t \in U}|t|) 
\end{array}
\right.
\end{equation}
for appropriately chosen constants $c_{1}$ and $c_{2}$. 
By virtue of the conditions on $a,b,h_{-},h_{+}$ it follows that
$\sup_{t \in U}\| k_{1}(t,\cdot)\|_{a,b}< \infty$. 
Since $X(\cdot,u) = c(\alpha(u))^{-1}L_{\alpha(u),h(u)}(\cdot)$,
Theorem 9.2(a) gives $h(u)$-localisability of $Y$ with
$Y_{u}'(\cdot) = X_{u}'(\cdot,u) = c(\alpha(u))^{-1}(L_{\alpha(u),h(u)})'_{u}(\cdot)
=  c(\alpha(u))^{-1}L_{\alpha(u),h(u)}(\cdot)$.

For part (b), we choose $U$ and the numbers $a,b,h_{-},h_{+}$ 
to satisfy the conditions stipulated in the proof of (a) but also to
satisfy $h_{-} >1/a +(1/a - 1/b) >0$,  so in particular 
$h(v)-1/\alpha(v') > 0$ for all $v,v' \in U$.
Again as in Theorem 7.4, 
$$
|f(t,v,x)-f(t',v,x)|,
|f_{v}(t,v,x)
-f_{v}(t',v,x)| 
\leq  k_{2}(t,t',x)
$$
for $t, v \in U, x\in \bbbr$, where
\begin{equation}
k_{2}(t,t',x) = \left\{
\begin{array}{ll}
    c_{3}|t-t'|^{h_{-}-1/a}
    & ( |x-\textstyle{\frac{1}{2}}(t-t')| \leq |t-t'|) \\
    c_{4}|x-\textstyle{\frac{1}{2}}(t-t')|^{h_{+}-1/b-1}|t-t'| 
     & ( |x-\textstyle{\frac{1}{2}}(t-t')| > |t-t'|) 
\end{array} \right.
\end{equation}
for constants $c_{3},c_{4}$.  Then 
$\| k_{2}(t,t',\cdot)\|_{a,b}\leq 
c_{5}|t-t'|^{1/a}$.
The conditions of Theorem 9.2(c) are 
satisfied with $\eta = 1/a>1/\alpha(u)$, so strong localisability follows.
\Box
\medskip

Note that the differentiability conditions in Theorem 9.3 could be
weakened slightly to H\"{o}lder conditions for which Theorem 9.2
would still be applicable.

Recall that an $\alpha$-stable L\'{e}vy motion, $0<\alpha<2$, 
is a process of $D(\bbbr)$ with
stationary independent increments which have a strictly $\alpha$-stable
distribution.   Taking $M$ as a symmetric $\alpha$-stable 
random measure on $\bbbr$,  the $\alpha$-stable L\'{e}vy 
motion may be represented as
\begin{equation} 
L_{\alpha}(t)= M[0,t]=\int \1_{[0,t]}(x)M(dx) = c(\alpha)\sum_{(\X,\Y)\in
\Pi} \1_{[0,t]}(\X) \Y^{<-1/\alpha>}\quad (t \in \bbbr),
\label{aslm2}
\end{equation}
where $\Pi$ is the Poisson process on $\bbbr^{2}$ with ${\cal L}^{2}$ 
as mean measure, $\1_{[0,t]}$ is the indicator function and  $c(\alpha) =
\left(2\alpha^{-1}\Gamma(1-\alpha)
\cos(\textstyle{\frac{1}{2}}\pi
\alpha)\right)^{-1/\alpha}$. 
Then $L_{\alpha}$ is $1/\alpha$-sssi and is strongly $1/\alpha$-localisable
in $D(\bbbr)$.

\begin{theo} {\it (Multistable L\'{e}vy motion)}.
Let $\alpha: \bbbr \to (0,2)$  and $a: \bbbr \to \bbbr^{+}$ 
be continuously differentiable, and 
define 
\begin{equation} 
Y(t)  = a(t)c(\alpha(t))\sum_{(\X,\Y)\in\Pi}
\1_{[0,t]}(\X)\Y^{<-1/\alpha(t)>} \quad (t \in \bbbr).
\label{aslm3}
\end{equation}

(a) If $1<\alpha(u) <2$ then  $Y$ is  $1/\alpha(u) $-localisable at
$u$, 
with $Y_{u}' = a(u)L_{\alpha(u)}$.

(b) If $0<\alpha(u) <1$ and $\alpha'(u) \neq 0$
then  $Y$ is  $1$-localisable at
$u$ 
with $\{Y_{u}'(t): t \in \bbbr\} = \{tW: t \in \bbbr\}$,
where $W$ is the random variable 
$$W = a(u)c(\alpha(u))
\sum_{(\X,\Y)\in\Pi}
\1_{[0,u]}(\X)\Y^{<-1/\alpha(u)>}\left(\frac{\alpha'(u)}{\alpha(u)^{2}}|\log|\Y||
+\frac{d}{du}(a(u)c(\alpha(u)))\right).$$
\end{theo}

\noindent{\it Proof.}  (a) By Proposition 3.2 
the term $a(t)c(\alpha(t))$ does not affect localisability. 
Define a random field by
$$
 X(t,v) =  \sum_{(\X,\Y)\in\Pi} \1_{[0,t]}(\X) 
 \Y^{<-1/\alpha(v)>} \quad (t,v \in \bbbr).
$$
Taking  $f(t,v,x) =\1_{[0,t]}(x)$ the conditions of Theorem 9.2(a)
are satisfied with $h = 1/\alpha(u) < 1\equiv \eta$, 
so the result follows from the localisability of 
$L_{\alpha}$. 

(b) In the case where $a(t)c(\alpha(t)) =1$
\begin{eqnarray*} \frac{Y(u+rt)-Y(u)}{r} &=&
\frac{1}{r}\sum_{(\X,\Y)\in\Pi}
\1_{[0,u]}(\X)(\Y^{<-1/\alpha(u+rt)>}
-\Y^{<-1/\alpha(u)>})\\
&&+\frac{1}{r}\sum_{(\X,\Y)\in\Pi}
\1_{[u,u+rt]}(\X)\Y^{<-1/\alpha(u+rt)>}.
\end{eqnarray*} 
Letting $r \to 0$ the second term vanishes if $1/\alpha(u) >1$ and the
first term converges to $W$ in finite dimensional distributions.
The general case is similar.
\Box 

\medskip
Note that Theorem 9.4(b) illustrates a general phenomenon that occurs
when the process $\{X(t,u): t \in \bbbr\}$ is 
$h(u)$-localisable at $u$ where $h(u)>1$. The process
$\{Y_{u}'(t): t \in \bbbr\}$ will typically be 
$1$-localisable at $u$, with the dominant
component of $Y_{u}'(t)$
derived from $(X(u+rt,u+rt)-X(u+rt,u))/r$ rather than from 
$X_{u}'(t,u)$.

As explained in \cite[Section 7.6]{ST}, there are two ways to extend the linear 
fractional stable motion to the
case $H = 1/ \alpha$. Apart from the L\'evy motion considered above, 
one may define the following process, called log-fractional stable motion:
\begin{equation}\label{lfsm}
\Lambda_\alpha(t)= \int_{-\infty}^{\infty} \left( \log(|t-x|) - \log(|x|)\right)
M(dx) \quad (t \in \bbbr)
\end{equation}
where, as usual, $M$ is an $\alpha$-stable random measure. 
This process is well-defined only for $\alpha \in (1,2]$ 
(the integrand does not belong to ${\cal F}_{\alpha}$ for $\alpha \leq 1$). 
It is $1/\alpha$-self-similar with stationary increments. 
Unlike the L\'evy motion, however, its increments are not independent.
Another difference is that 
log-fractional stable motion does not have a version in $D(\bbbr)$,
so we cannot speak of strong localisability. 

\begin{theo} (Log-fractional multistable motion)
Let $\alpha : \bbbr \to (1, 2)$ and $a$ be continuously differentiable, and define
\begin{equation}\label{lfmsm}
Y(t)= a(t) \sum_{(\X,\Y) \in \Pi} \left( \log|t-\X| -
\log|\X|\right) \Y^{<-1/\alpha(t)>} \quad (t \in \bbbr).
\end{equation}
Then $Y$ is $1/\alpha(u)$-localisable at all $u \in \bbbr$, 
with $Y'_u = a(u) \Lambda_{\alpha(u)}$.
\end{theo}

\noindent{\it Proof.} 
The proof is similar to that of Theorem 9.4, 
by considering the field
\begin{equation}
X(t,v)= \sum_{(\X,\Y) \in \Pi} \left( \log|t-\X| - \log|\X|\right) 
\Y^{<-1/\alpha(v)>}  \quad (t,v \in \bbbr),
\end{equation} 
with Theorem 9.2(a) is applied to $f(t,v,x)= \log|t-x| - \log|x|$.
\Box
\medskip

For a final example we give a multistable version of Theorem 7.6

\begin{theo} (Multistable reverse Ornstein-Uhlenbeeck
process)
Let $\lambda>0$ and $\alpha : \bbbr \to (1, 2)$ be continuously
differentiable.  Let
\begin{equation*}
Y(t)= \sum_{(\X,\Y) \in \Pi, \X \geq t} 
\exp(-\lambda(\X-t)) \Y^{<-1/\alpha(t)>}  \quad (t \in \bbbr).
\end{equation*}
Then $Y$ is $1/\alpha(u)$-localisable at all $u \in \bbbr$, 
with $Y'_u = c(\alpha(u))^{-1}L_{\alpha(u)}$, where $L_{\alpha}$ is L\'{e}vy
$\alpha$-stable motion.  
\end{theo}

\noindent{\it Proof.} 
Let
\begin{equation}
X(t,v) = \sum_{(\X,\Y) \in \Pi, \X \geq t} \exp(-\lambda(\X-t))
\Y^{<-1/\alpha(v)>}  \quad (t,v \in \bbbr).
\end{equation} 
Then for each $v$, $X(\cdot,v)$ is everywhere $1/\alpha(v)$-localisable 
with $X_{u}'(\cdot,v) = c(\alpha(v))^{-1}L_{\alpha(v)}(\cdot)$ by Theorem 7.6.  
Applying Theorem 9.2(a) with
$f(t,v,x)= \1_{[t,\infty)}\exp(-\lambda(x-t))$
gives the conclusion.
\Box


\section{Further work}
\setcounter{equation}{0}
\setcounter{theo}{0}
\medskip

There are a great many possible variants and extensions of this
work. Localisable processes of many other forms may be constructed. For example
multistable processes with skewness and the class of
 stationary localisable processes deserve investigation.  There may be 
 advantages in seeking other representations of multistable
processes  such as by sums involving
arrival times of a Poisson process or as stochastic integrals with
respect to suitably constructed random measures. 
Our conditions for localisability could certainly be weakened and further 
techniques for establishing localisability and in particular strong
localisability developed.  It would be interesting to study long 
range dependence of Multistable processes. Effective techniques for simulation 
and inference on parameters for these processes are also needed.  
We will be addressing some of these matters in a sequel to this paper.

\end{document}